\documentclass[12pt]{article}

\usepackage{amsthm}

\usepackage{epsfig}
\usepackage{latexsym}
\usepackage{amsmath}
\usepackage{amssymb}
\usepackage{amsfonts}
\usepackage[size=1.5em,nohug,heads=LaTeX,PostScript=dvips]{diagrams}

\usepackage{epstopdf}

\def\1{{\mathbf 1}}

\theoremstyle{definition}

\newtheorem{dfn}{Definition}[section]
\newtheorem{defn}[dfn]{Definition}
\newtheorem{definition}[dfn]{Definition}

\newtheorem{rem}[dfn]{Remark}
\newtheorem{remark}[dfn]{Remark}

\theoremstyle{plain}

\newtheorem{thm}[dfn]{Theorem}
\newtheorem{theorem}[dfn]{Theorem}

\newtheorem{lem}[dfn]{Lemma}
\newtheorem{lemma}[dfn]{Lemma}

\newtheorem{prop}[dfn]{Proposition}
\newtheorem{proposition}[dfn]{Proposition}

\newtheorem{problem}[dfn]{Problem}

\newtheorem{cor}[dfn]{Corollary}
\newtheorem{corollary}[dfn]{Corollary}

\newtheorem{notation}[dfn]{Notation}

\def\k{\mathbf k}

\def\proof{\par\medskip\noindent{\it Proof. }}

\def\C{{\mathbb C}}
\def\R{{\mathbb R}}

\def\Z{{\mathbb Z}}
\newcommand{\K}{{\mathcal K}}

\def\P{{\mathbb P}}

\def\N{{\mathbb N}}

\def\al{\alpha}
\def\be{\beta}

\def\del{\delta}
\def\De{\Delta}
\def\Si{\Sigma}
\def\si{\sigma}

\def\la{\lambda}
\def\La{\Lambda}

\def\acts{\curvearrowright}

\def\embed{\hookrightarrow}
\def\8{\infty}
\def\<{\langle}
\def\>{\rangle}

\def\ov{\overrightarrow}
\def\ol{\overline}

\begin{document}

\title{Stability inequalities and universal Schubert calculus of rank $2$}
\author{Arkady Berenstein and Michael Kapovich}
\date{August 10, 2010}

\maketitle

\begin{abstract} The goal of the paper is to introduce a version of Schubert calculus for each dihedral reflection group $W$.
That is,  to each ``sufficiently rich'' spherical building $Y$ of type $W$ we associate a certain cohomology 
theory $H^*(Y)$ and verify that, first, it  depends only on $W$
(i.e., all such buildings are ``homotopy equivalent'') and second, $H^*(Y)$ 
is the associated graded of the coinvariant algebra of $W$ under certain filtration. 
We also construct the dual homology ``pre-ring'' on $Y$. The convex ``stability'' cones in $(\R^2)^m$ defined 
via these (co)homology theories of $Y$ are then shown to solve the problem of classifying weighted semistable 
$m$-tuples on $Y$ in the sense of \cite{KLM1}; equivalently, they are cut out by the {\em generalized triangle inequalities} for 
thick Euclidean buildings with the Tits boundary $Y$.  
The independence of the (co)homology theory of $Y$ refines the result of \cite{KLM2}, which asserted that the 
stability cone depends on $W$ rather than on $Y$. Quite remarkably, the cohomology ring 
$H^*(Y)$ is obtained from a certain universal algebra $A_t$ by a kind of ``crystal limit'' 
that has been previously introduced by Belkale-Kumar for  the cohomology of flag varieties and Grassmannians. 
Another degeneration of  $A_t$ leads to the homology theory $H_*(Y)$. 
\end{abstract}

\section{Introduction}

Alexander Klyachko in \cite{Kl} solved the old problem on
eigenvalues of sums of hermitian matrices: His solution was to
interpret the eigenvalue problem as an existence problem for
certain parabolically stable bundles over $\C \P^1$, so that the
inequalities on the eigenvalues are stated in terms of the
Schubert calculus on Grassmannians. Klyachko's work was later
generalized by various authors to cover general semisimple groups,
see e.g. \cite{BS, KLM1}. The stable bundles were replaced in
\cite{KLM1}  with {\em semistable weighted configurations} on
certain spherical buildings and  the eigenvalue problem was
interpreted as a {\em triangle inequalities problem} for
the vector-valued distance function on nonpositively curved symmetric
spaces and Euclidean buildings. Still, the solution depended
heavily (and was formulated in terms of) Schubert calculus on
generalized Grassmannians $G/P$, where $G$ is a (complex or
real) semisimple Lie group and $P$'s are maximal parabolic
subgroups of $G$.

The present work is a part of our attempt to generalize Lie theory
to the case of nonexistent Lie groups having non-crystallographic
dihedral groups $W=I_2(n)$ (of order $2n$) as their Weyl groups. For such Weyl
groups, one cannot define $G$ and $P$, but there are
spherical (Tits) buildings $Y$, whose vertex sets serve as generalized
Grassmannians $G/P$. Moreover, we also have thick discrete
and nondiscrete Euclidean buildings for the groups $I_2(n)$ (see
\cite{BK}), so both problems of existence of semistable weighted
configurations and computation of triangle inequalities for
vector-valued distance functions on Euclidean buildings certainly
make sense. The goal of this paper is to compute these
inequalities (by analogy with \cite{BS, KLM1}) in terms of the
Borel model for $H^*(G/P)$ and to verify that they solve the
problem of existence of semistable weighted configurations and the
equivalent problem of computation of triangle inequalities in the
associated affine buildings.

Our main results can be summarized as follows:

Let ${\mathfrak Y}$ be a rank 2 affine building with the Weyl group $W=I_2(n)$, let $\De$ denote the
positive Weyl chamber for $W$. We then obtain a $\De$-valued distance function $d_\De(x,y)$
between points $x, y\in {\mathfrak Y}$, see \cite{KLM1} or \cite{KLM3}. Then

\begin{thm}\label{main0}
There exists a geodesic $m$-gon $x_1\cdots x_m$ in ${\mathfrak Y}$ with the
$\De$-side-lengths $\la_1,...,\la_m$ if and only if the vectors
$\la_1,...,\la_m$  satisfy the {\em Weak Triangle Inequalities}
(the {\em stability inequalities}):
\begin{equation}\label{wti}
w(\la_i - \la_j^*) \le_{\De^*} \sum_{k\ne i, k\ne j} \la_k^*,
\quad w\in W
\end{equation}
taken over all distinct $i, j\in \{1,...,m\}$.
\end{thm}

Here $\la^*=-w_\circ(\la)$ 
is the vector contragredient to $\la$
($w_\circ\in W$ is the longest element). The order $\le_{\De^*}$ is
defined with respect to the obtuse cone $\De^*$ dual to $\De$:
$$
\De^*=\{\nu: \nu\cdot \la\ge 0, \forall \la\in \De\}.
$$
(Recall that $\la \le _{\De^*} \nu \iff \nu -\la \in \De^*$.)

The key idea behind the proof is that although we do not have
smooth homogeneous manifolds $G/P$, we still can define some
kind of Schubert calculus on the sets of ``points'' $Y_1$ and
``lines'' $Y_2$ in appropriately chosen Tits buildings $Y$ (replacing
$G/P$'s). We define certain ``homology pre-rings'' $H_*(Y_l, \widehat\k)$,  
$l=1, 2$, (``Schubert pre-calculus'') which reflect the 
intersection properties of ``Schubert cycles'' in $Y_l$.
We then show that this calculus is robust enough to solve the
existence problem for weighted semistable configurations.

\medskip
We then promote the cohomology pre-rings to rings. To this end, we introduce
the {\em universal Schubert calculus}, i.e., we define a {\em cohomology ring} $H^*(Y, \k)=A_t$ 
for each  reflection group of rank 2, based on a generalization of the Borel model  
for the computation of cohomology rings of flag varieties. 
One novelty here is that in the definition of $A_t$ we allow $t\in \C^\times$, thereby providing an 
interpolation between cohomology rings of complex flag manifolds; for $t$ a primitive $n$-th root of unity,  
$A_t$ defines  $H^*(Y, \k)$, the cohomology rings of the buildings $Y$ with the Weyl group $W=I_2(n)$. 
We, therefore, think of the family of rings $A_t$ as ``universal Schubert calculus'' in rank 2. 
An odd feature of the rings $A_t$ is that even for the values of $t$ which are roots of unity, the structure constants 
of $A_t$ are typically irrational ($t$-binomials), so we do not have a natural geometric model for these $A_t$. 
In order to link $A_t$ to geometry, we define a (trivial) deformation $A_{t,\tau}, \tau\in \R_+$, of $A_t$.   
Sending $\tau$ to $0$ we obtain an analogue of the Belkale-Kumar degeneration $H^*_{BK}(Y, \k)=gr(A_t)$
of $A_t$. On the other hand, by sending $\tau$ to $\infty$, we recover the pre-ring $H_*(Y, \widehat{\k})$ 
given by the  Schubert pre-calculus. Therefore, $A_{t}$ interpolates between  $H^*_{BK}(Y, \k)$ and 
$H_*(Y, \widehat{\k})$. The same relation holds for the cohomology rings of ``Grassmannians,'' 
$B^{(l)}_t=H^*(Y_l, \k)\subset A_t$, their Belkale-Kumar degenerations $H^*_{BK}(Y_l, \k)=gr(B^{(l)}_t)$ 
and pre-rings  $H_*(Y_l, \widehat{\k})$.  We then observe (Section \ref{sec:sti}) that the system  of 
strong triangle inequalities defined by $H_*(Y, \widehat{\k})$ also 
determines the stability cone $\K_m(Y)$ for the building $Y$.  
In \S \ref{generalize} we introduce systems of linear inequalities determined by certain based rings $A$, 
generalizing $A_t$. Specializing these inequalities to the case $A=A_t$, using the results of \S \ref{sec:sti}, 
we recover the stability cones $\K_m(Y)$. Therefore, the systems of 
inequalities defined by $A_t, B_t^{(l)}, gr(A_t), gr(B_t^{(l)})$ and $H_*(X)$ are all equivalent. 
In this section, we also prove that the system of Weak Triangle Inequalities, 
determined by $H^*_{BK}(Y_l, \k)$, equivalently, $H_*(Y_l, \widehat{\k})$, ($l=1,2$) is irredundant. 
This is remeniscent of the result by Ressayre who proved irredundancy of the Beklale-Kumar inequalities 
in the context of complex reductive groups.

\medskip
{\bf Acknowledgments.} Our collaboration on this project started
at the AIM workshop ``Buildings and Combinatorial Representation
Theory'' in 2007 and we are grateful to AIM for this opportunity.
The first author was supported by the NSF grant DMS-08-00247. The
second author was supported by the NSF grants DMS-05-54349 and DMS-09-05802.

\tableofcontents

\section{Coxeter complexes}

Let $A$, the {\em apartment}, be either the Euclidean space $E=E^N$
or the unit $N-1$-sphere $S=S^{N-1}\subset E^N$ (we will be
primarily interested in the case of Euclidean plane and the
circle). If $A=S$, a {\em Coxeter group} acting on $A$ is a
finite group $W$ generated by isometric reflections. If $A=E$, a
{\em Coxeter group} acting on $A$ is a group $W_{af}$ generated by
isometric reflections in hyperplanes in $A$, so that the linear
part of $W_{af}$ is a Coxeter group acting on $S$. Thus,
$W_{af}=\La \rtimes W$, where $\La$ is a certain (countable or
uncountable) group of translations in $E$. We will use the notation $\1$ for the identity in $W$ and $w_\circ$ for the longest element of $W$ 
with respect to the word-length function $\ell: W\to \Z$ with respect to the standard Coxeter generators $s_i$.

\begin{defn}
A spherical or Euclidean\footnote{Also called {\em affine}.}
{\em Coxeter complex} is a pair $(A, G)$, of the form $(S, W)$ or $(E, W_{af})$. 
The number $N$ is called the {\em rank} of the Coxeter complex.
\end{defn}

A spherical Coxeter complex $(S,W)$ is {\em essential}  
if $W$ has no global fixed points in $S$.

A  {\em wall} in the Coxeter complex $(A, G)$ is the fixed-point
set of a reflection in $G$. A {\em half-apartment} in $A$ is a
closed half-space bounded by a wall. A {\em regular point} in a
Coxeter complex is a point which does not belong to any wall. A
{\em singular point} is a point which is not regular.

\begin{rem}
Note that in  the spherical case, there is a natural cell complex
in $S$ associated with $W$. However, the affine case, when
$W_{af}$ is nondiscrete, there will be no natural  cell complex 
attached to $W_{af}$.
\end{rem}

{\em Chambers} in $(S, W)$ are the fundamental domains for the
action $W\acts S$, i.e., the closures of the connected components
of the complement to the union of walls. We will use the notation 
$\De_{sph}$ for a fixed (positive) fundamental domain.  

An {\em affine} Weyl chamber in $(A, W_{af})$ is a fundamental domain $\De=\De_{af}$ 
for a conjugate $W'$ of $W$ in $W_{af}$, i.e. it is a cone over $\De_{sph}$ with the tip at a point $o$ fixed by $W'$.  

A {\em vertex} in $(A, G)$ is a (component of, in the spherical
case) the 0-dimensional intersection of walls. We will consider
almost exclusively only those Coxeter complexes which have at
least one vertex; such complexes are called {\em essential}.
Equivalently, these are spherical complexes where the group $G$
does not have a global fixed point and those Euclidean Coxeter
complexes where $W$ does not have a fixed point in $S$.

In the spherical case, the notion of {\em type} is given by the projection
$$
\theta: S\to S/W=\De_{sph},
$$
where the quotient is the spherical Weyl chamber.

Let $s_i\in W$ be one of the Coxeter generators. We define the
{\em relative length functions} $\ell_i$ on $W$ as follows: 
$\ell_i(w)$ is the length of the shortest element of the coset $w\<s_i\>\subset W /\<s_i\>$.
In the case when $W$ is a finite dihedral group, $\ell_i(w)$ equals is the combinatorial 
distance from the vertex $w(\zeta_i)$
to the positive chamber $\Delta_{sph}$ in the spherical Coxeter complex $(S^1, W)$.
Here, $\zeta_i$ is the vertex of $\Delta_{sph}$ fixed by $s_i$.

\section{Metric concepts}\label{metric}

\begin{notation}
Let $Y, Z$ be subsets in a metric space $X$. 
Define the (lower) distance $d(Y, Z)$ as 
$$
\inf_{y\in Y, z\in Z} d(y, z). 
$$
If $Z$ is a singleton $\{z\}$, we abbreviate $d(\{z\}, Y)$ to $d(z, Y)$. 
In the examples we are interested in, the above infimum is always realized.

For a subset $Y\subset X$, we let 
$B_r(Y)$ denote the {\em closed $r$-neighborhood} of $Y$ in $X$,
i.e.,
$$
B_r(Y):= \{x\in X: d(x, Y)\le r 
\}.
$$
For instance, if $Y=\{y\}$ is a single point, then $B_r(Y)=B_r(y)$
is the closed $r$-ball centered at $y$. Similarly, we define ``spheres centered at $Y$'' 
$$
S_r(Y):= \{ x\in X: d(x, Y)=r\}. 
$$ 
\end{notation}

A metric space $X$ is called {\em geodesic} if  every two points
in $X$ are connected by a (globally distance-minimizing) geodesic.
Most metric spaces considered in this paper will be geodesic.
Occasionally, we will have to deal with metrics on disconnected
graphs: In this case we declare the distance between points in
distinct connected components to be infinite.

For a pair of points $x, y$ in a metric space $X$  we let
$\ol{xy}$ denote a closed geodesic segment (if it exists) in $X$
connecting $x$ and $y$. As, most of the time, we will deal with
spaces where every pair of points is connected by the unique
geodesic, this is a reasonable notation.

We refer to \cite{Ballmann} or \cite{BH} for the definition of a CAT($k$) metric space. 
We will think of the distances in CAT(1) spaces as {\em angles} 
and, in many cases, denote these distances $\angle(xy)$.

\medskip
The following characterization of 1-dimensional CAT(1) spaces will
be important:

A 1-dimensional metric space (a metric graph) is a CAT(1) space if
and only if the length of the shortest embedded circle in $X$ is
$\ge 2\pi$.

If $G$ is a metric graph, where each edge is given the length
$\pi/n$, then the CAT(1) condition is equivalent to the assumption
that girth of $G$ is $\ge 2n$.

Fix an integer $n\ge 2$. Similarly to \cite{BK}, a type-preserving map of
bipartite graphs $f: G\to G'$ is said to be $(n-1)$-isometric if:

1. $\forall x, y\in V(G), d(x,y)< n-1 \Rightarrow d(f(x),
f(y))=d(x,y)$.

2. $\forall x, y\in X, d(x,y)\ge n-1 \Rightarrow d(f(x), f(y))\ge
n-1$.

Here $d$ is the combinatorial path-metric on $G$, which is
allowed to take infinite values on points which belong to distinct
connected components. One can easily verify that the concept of an
$(n-1)$-isometric map is equivalent to the notion of a 
type-preserving map graphs which preserves the {\em bounded
distance} on the graphs defined in \cite{Tent}.

\section{Buildings}

{\bf Spaces modeled on Coxeter complexes.}

Let $(A,G)$ be a Coxeter complex (Euclidean or spherical).

\begin{defn}\label{modeled}
A space  {\em modeled on the Coxeter complex $(A, G)$} is
a metric space $X$ together with an atlas where charts are
isometric embeddings $A\to X$ and the transition maps are
restrictions of the elements of $G$. The maps $A\to X$ and
their images are called {\em apartments} in $X$. Note that (unlike
in the definition of an atlas in a manifold) we do not require the
apartments to be open in $X$.
\end{defn}

Therefore, all $G$-invariant  notions 
defined in $A$, extend to $X$. In particular, we will talk about
vertices, walls, chambers, etc. 

\begin{notation}
We will use the notation $\De_i$ for chambers in spherical buildings. 
\end{notation}

{\em Rank} of $X$ is the rank of the
corresponding Coxeter complex. 


A space $X$ modeled on $(A,G)$ is called {\em discrete} if the group $G$ is discrete.
This is automatic in the case of spherical Coxeter complexes since $G$ is finite in this case.

\medskip
A {\em spherical building} modeled on $(S, W)$ is a CAT(1) space $Y$ modeled on
$(S, W)$ which satisfies the following condition:

{\bf Axiom (``Connectedness'').} Every two points $y_1, y_2\in Y$ are contained in a common apartment.

\medskip
The group $W$ is called the {\em Weyl group} of the spherical
building $Y$.

Spherical buildings of rank 2 (with the Weyl group of order $\ge
4$) are called {\em generalized polygons}. They can be described
combinatorially as follows:

A building $Y$ is a bipartite graph of girth $2n$ and valence $\ge 2$ at every vertex,
so that every two vertices are connected by a path of the combinatorial length $\le n$.
To define a metric on $Y$, we identify each edge of the graph with the segment of length $\pi/n$.

\medskip
A {\em Euclidean (or, affine) building} modeled on $(A, W_{af})$ is a CAT(0) space $X$ modeled on
$(A, W_{af})$ which satisfies the following conditions:

{\bf Axiom 1.}
(``Connectedness'') Every two points $x_1, x_2\in X$ belong to a common apartment.

\medskip
{\bf Axiom 2.} There is an extra axiom (comparing to the spherical
buildings) of ``Angle rigidity'', which will be irrelevant for the
purposes of this paper. It says that for every $x\in X$, the space
of directions $Y=\Si_x(X)$ satisfies the following:
$$
\forall \xi, \eta\in Y, \quad \angle(\xi, \eta)\in W\cdot \angle(\theta(\xi), \theta(\eta)).
$$
Here $\theta: Y\to \De_{sph}$ is the {\em type projection}. We refer to \cite{BK, KL, Parreau} for the details.
Note that Axiom 2 is redundant in the case of discrete Euclidean buildings.

The finite Coxeter group $W$ (the linear part of $W_{af}$) is called the {\em Weyl group} of the Euclidean building $X$.

\medskip
A building $X$ is called {\em thick} if every wall in $X$ is the intersection of
(at least) three half-apartments.

\bigskip
We now specialize our discussion of  buildings to the
case of rank 2 (equivalently, $1$-dimensional) spherical buildings.

Chambers $\De_1, \De_2$ in a spherical building $Y$ are called
{\em antipodal} if the following holds. Let $A\subset Y$ be an 
apartment containing both $\De_1, \De_2$ (it exists by the
Connectedness Axiom). Then $\De_1=-\De_2$ inside $A$. More generally, if $Y$ is a bipartite graph of diameter $n$, 
then two edges $e_1, e_2$ of $Y$ are called {\em antipodal} if the minimal distance between vertices of these edges is exactly $n-1$.

\medskip
Let $W=I_2(n)$ be the dihedral group of order $2n$. We regard {\em
type}  of a vertex $x$ (denoted $type(x)$) of a bipartite graph
(in particular, of a spherical building with the Weyl group $W$)
to be an element of $\Z/2$. We let $W_l, l=1, 2$ denote the
stabilizer of the vertex of type $l$ in the positive (spherical)
chamber $\De_+$ of $W$.

Let $Y$ be a rank $2$ spherical building with the Weyl group
$W$; we will use two metrics on $X$:

1. The combinatorial path-metric $d=d_Y$ between the vertices of
$X$, where each edge is given the unit length. This metric extends
naturally to the rest of $Y$: we will occasionally use this fact.

2. The (angular) path metric $\angle$ on $Y$ where every edge has
the length $\pi/n$. Given a subset $Z\subset Y$ we let $B_r(Z)$
and $S_r(Z)$ denote the closed $r$-ball and $r$-sphere in $X$ with
respect to the combinatorial metric.

 The building $Y$ has two vertex types identified with $l\in \Z/2$; accordingly, the vertex set of $Y$ is
 the disjoint union $Y_0\cup Y_1$ of the {\em Grassmannians} of type $l=0, 1$. When $l$ is fixed, by abusing the
 notation, we will denote by $B_r(Z)$ (and $S_r(Z)$) the intersection of the corresponding ball
 (or the sphere) with $Y_l$.  We will only use these concepts when $Z$ is a vertex or a
chamber $\De$ of a spherical building. The balls $B_r(\De)\subset
Y_l$ will serve as {\em Schubert cycles} in the Grassmannian
$Y_l$, while the spheres $S_r(\De)$ will play the role of (open)
Schubert cells.

\section{Weighted configurations and geodesic polygons}

{\bf Weighted configurations.}
Let $Y$ be a spherical building modeled on $(S, W)$. We recall that $\angle$ denotes the metric on $Y$. 
Given a collection $\mu_1,...,\mu_m$ of
non-negative real numbers (``weights'') we define a {\em weighted configuration} on $Y$ as a map
$$
\psi: \{1,...,m\}\to Y,\quad \psi(i)=\xi_i\in Y.
$$
We thus get $n$ points $\xi_i, i=1,...,m$ on $Y$ assigned the weights $\mu_i, i=1,...,m$.
We will use the notation
$$
\psi=(\mu_1 \xi_1,...,\mu_m\xi_m).
$$

Let $\theta: Y\to \De_{sph}$ denote the {\em type-projection}
to the spherical Weyl chamber.
Given a weighted configuration $\psi=(\mu_1 \xi_1,...,\mu_m\xi_m)$ on $Y$, we define
$\theta(\psi)$, the {\em type of $\psi$} to be the $n$-tuple of vectors
$$
(\la_1,...,\la_m)\in \De^m,
$$
where $\la_i=\mu_i \theta(\xi_i)$.

Following \cite{KLM1}, for a finite weighted configuration
$\psi=(\mu_1 \xi_1,...,\mu_m\xi_m)$ on $Y$, we define the function
$$
slope_\psi: Y\to \R,
$$
$$
slope_\psi(\eta)=-\sum_{i=1}^m \mu_i \cos(\angle(\eta, \xi_i)).
$$

\begin{defn}
A weighted configuration $\psi$ is called {\em semistable} if the associated slope function is $\ge 0$ on $Y$.
\end{defn}

It is shown in \cite{KLM1} that $slope_\psi$ coincides with the
Mumford's numerical stability function for weighted configurations
on generalized flag-varieties. What's important is the fact that
the above notion of stability, unlike the stability conditions in
algebraic and symplectic geometry, does not require a group action
on a smooth manifold. (Actually, it does not need any group at
all.)

\medskip
{\bf Vector-valued distance functions.} Let $X$ be a Euclidean
building modeled on $(A, W_{af})$. Our goal is to define the
$\De$-valued distance function $d_\De$ on $X$, where $\De=\De_{af}$ is an (affine) Weyl chamber of $(A, W_{af})$. 
We first define this function on $A$. Let $o\in A$ denote the point fixed by $W$.
We regard $o$ as the origin in the affine space $A$, thus giving
$A$ the structure of a vector space $V$. Then, given two points
$x, y\in A$, we consider the vector $v=\ov{xy}$ and project it to
a vector $\bar{v}\in \De$ via the map
$$
V\to V/W=\De.
$$
Then $d_\De(x,y)=\bar{v}$. It is clear from the construction that $W_{af}$ preserves $d_\De$. 
Suppose now that $x, y\in X$. Then there exists an apartment
$\phi: A\to X$ whose image contains $x$ and $y$. We then set
$$
d_\De(x,y):= d_\De(\phi^{-1}(x), \phi^{-1}(y))\in \De.
$$
Since the transition maps between the charts are in $W_{af}$, it follows that the distance function
$d_\De$ on $X$ is well-defined. Note that $d_\De$ is, in general, non-symmetric:
$$
d_\De(x,y)=\la \iff d_\De(y,x)=\la^*, \quad \la^*=-w_\circ(\la),
$$
where $w_\circ\in W$ is the longest element. Hence, unless $w_\circ=-1$, $d_\De(x,y)\ne d_\De(y,x)$.

A {\em (closed) geodesic $m$-gon} on $X$ is an $m$-tuple of points
$x_1,...,x_m$, the {\em vertices} of the polygon. Since for every
two points $x, y\in X$ there exists a unique geodesic segment
$\ol{xy}$ connecting $x$ to $y$, the choice of vertices uniquely
determines a closed 1-cycle in $X$, called {\em geodesic polygon}.
We will use the notation $x_1\cdots x_m$ for this polygon. The
{\em $\De$-side-lengths} of this polygon are the vectors
$\la_i=d_{\De}(x_i, x_{i+1})$, where $i$ is taken modulo $m$.

\medskip
The following  is proven in \cite{KLM2}:

\begin{thm}\label{T1}
Let $Y$ be a thick spherical building modeled on $(S, W)$ and $X$ be a thick Euclidean building modeled on
$(A, W_{af}=\La \rtimes W)$, for an arbitrary $\La$. Then:

There exists a weighted semistable configuration $\psi$ of type $(\la_1,...,\la_m)$ on $Y$ if and only if
there exists a closed geodesic $m$-gon $x_1...x_m$ in $X$ with the $\De$-side-lengths $(\la_1,...,\la_m)$.
\end{thm}

In particular, the existence of a semistable configuration (or a 
geodesic polygon) depends only on $W$ and nothing else. The way it
will be used in our paper  is to construct special spherical
buildings modeled on $(S^1, I_2(n))$ (buildings satisfying {\bf Axiom
A}), to which certain ``transversality arguments'' from
\cite{KLM1} apply.

\begin{definition}
Given a thick spherical building $X$ with the Weyl group $W$, we let $\K_m(X)$ denote the set of vectors 
$\ov\la=(\la_1,...,\la_m)$ in $\De^m$, so that  $X$ contains  a semistable weighted configuration of the type
$\ov\la$.  We will refer to $\K_m(X)$ as the {\em Stability Cone} of $X$. 
(These cones are also known as {\em Eigenvalue Cones} in the context of Lie groups and Lie algebras.) 
When $W$ is fixed, we will frequently abbreviate $\K_m(X)$ to $\K_m$ since this cone depends only on the dihedral group $W$. 
\end{definition}

Note that conicality  of $X$ is clear since a positive multiple of a semistable weighted configuration is again semistable. 
What is not obvious is that $\K_m(X)$ is a convex polyhedral cone. We will see in \S \ref{sec:sineq} 
(as a combination of the results of this paper and \cite{KLM1}) that this indeed always the case.

\section{$m$-pods}

Fix an integer $n\ge 2$. Let $r_1,...,r_m$ be positive integers
such that
$$
r_i+r_j\ge n, \forall i\ne j.
$$
Given this data, we define an $m$-pod $T$ as follows.

Let $B$ denote the bipartite graph which is the disjoint union of
the edges $\De_1$, $...$, $\De_m$. These edges will be the {\em bases}
of the $m$-pod $T$. Add to $B$ the vertex $z$ of type $l$, the
{\em center} of $T$.  Now, connect $z$ to the appropriate vertices
$x_i\in \De_i$ by the paths $p_i$ of the combinatorial lengths
$r_i$, so that
$$
r_i\equiv  type(x_i)+ type(z) \quad (\hbox{mod ~~2}),  i=1, ...,m.
$$
(The above equation uniquely determines each $x_i$.) The resulting graph is $T$.
The paths $p_i$ are the {\em legs} of $T$. It is easy to define the type of the vertices of $T$
(extending those of $x_1, ..., x_m, z$), so that $T$ is a bipartite graph.

\begin{figure}[tbh]
\centerline{\epsfxsize=4.5in \epsfbox{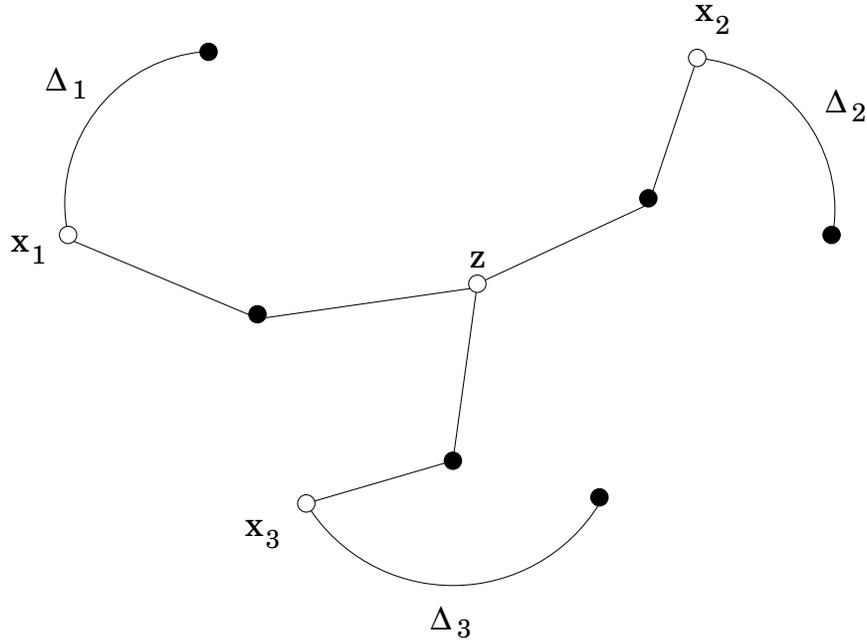}}
\caption{\sl Attaching a tripod.}
\label{F1}
\end{figure}

Suppose now that $Y$ is a bipartite graph of girth $\ge 2n$,
$\De_1,...,\De_m\subset Y$ are mutually antipodal edges and $r_1,...,r_m$
are positive integers so that
$$
r_i+r_j\ge n, \quad \forall i\ne j.
$$
We define a new bipartite graph $Y'$ by attaching the $m$-pod $T$
with legs of the lengths $r_i$, $i=1,...,m$, and with the bases
$\De_1,...,\De_m$.

\begin{lem}\label{L3}
The graph $Y'$ still has girth $\ge 2n$.
\end{lem}
\proof Since $r_i+r_j\ge n$ for $i\ne j$, the only thing we need
to avoid is having $i\ne j$ so that $r_i+r_j=n$ and $d_Y(x_i,
x_j)=n-1$ (since the edges $\De_k$ to which the $m$-pod $T$ is
attached are all antipodal). Suppose such $i, j$ exist. Then,
$$
type(x_i)+ type(x_j)\equiv r_i+r_j= n \quad (\hbox{mod~~} 2)
$$
and
$$
type(x_i)+ type(x_j)\equiv d_Y(x_i, x_j)=n-1 \quad (\hbox{mod~~} 2).
$$
Contradiction. \qed

\section{Buildings and free constructions}\label{sec:A}

We define a class of rank 2 spherical buildings $X$ with the
Weyl group $W=I_2(n)$ satisfying:

\medskip
{\bf  Axiom A.} 1. Each  vertex of $X$ has infinite valence, in particular, $X$ is thick. 

2. For each $m\ge 3$ the following holds. Let $\De_i,
i=1, ..., m$ be pairwise antipodal chambers in $X$ and let $0<r_i\le n-1$, 
$i=1,..., m$, be integers so that
$$
r_i+r_j\ge n, \quad \forall i\ne j.
$$
Then there exist infinitely many 
vertices $\eta\in X$ of both  types, so that
$$
d(\eta, \De_i)\le r_i.
$$
In other words, the intersection of metric spheres
$$
\bigcap_{i} S_{r_i}(\De_i)
$$
contains infinitely many vertices of both types.

\begin{rem}
1. For the purposes of the proof of Theorem \ref{main}, it suffices
to have this property for a fixed infinite collection $\De_1,
\De_2,...$ of pairwise antipodal chambers.

2. Clearly, Axiom fails for finite buildings. However, it also fails for some infinite buildings. For instance, 
it fails for the Tits buildings associated with the complex algebraic groups $Sp(4, \C)$ and $G_2(\C)$.  
\end{rem}

Buildings satisfying Axiom A constitute the class of ``sufficiently rich'' buildings mentioned in the Introduction: For these buildings we
will develop ``Schubert Calculus'' later in the paper.

\begin{lem}\label{anti}
Let $X$ be a thick rank 2 spherical building satisfying Axiom A. 
Let $\De_1,...,\De_m$ be pairwise antipodal chambers in $X$. Then there
exists a chamber $\De_{m+1}$ antipodal to all chambers $\De_1,...,\De_m$. 
\end{lem}
\proof Let $r_i:=n-1$. Then, by Axiom A, there exists a vertex $x\in X$ so that
$$
d(x, \De_i)=r_i, i=1,...,m. 
$$
For each $i$ we let $x_i\in \De_i$ be the vertex realizing $d(x, \De_i)$. Since $X$ has infinite valence at $x$, there exists a vertex $y\in
X$ incident to $x$, which does not belong to any of the geodesics $\ol{x x_i}, i=1,...,m$. It is then clear that $d(y, \De_i)=n-1$,
$i=1,...,m$. Therefore, the chamber $\De_{m+1}:=\ol{xy}$ is antipodal to all chambers $\De_1,...,\De_m$. \qed

\begin{rem}
It is not hard to prove that if $X$ is a thick building with Weyl group $I_2(n)$, then the conclusion of the 
above lemma holds for $m=2$ without any extra assumptions.  
\end{rem}

We now prove the existence of thick buildings satisfying Axiom A.

\begin{thm}\label{Aexists}
For each $n$ there exists a thick spherical building $X$ with Weyl group $W\cong I_2(n)$, countably many vertices 
and  satisfying Axiom A. Moreover, every 
(countable) graph of girth $\ge 2n$ embeds in a (countable) building satisfying Axiom A.   
\end{thm}
\proof 
We first recall the {\em free construction} of rank 2 spherical
buildings (see \cite{Tits, Ronan, FS}):

Let $Z$ be a connected bipartite graph of girth $\ge 2n$. Given
every pair of vertices $z, z'\in Z$ of distance $n+1$ from each
other, we add to $Z$ an edge-path $p$ of the combinatorial length  
$n-1$ connecting $z$ and $z'$; similarly, 
for every pair of vertices in $Z$ of distance $n$ from each other add an edge-path $q$  
of the combinatorial length $n$ connecting $z$ and $z'$. Let $\ol{Z}$ denote the graph 
obtained by attaching paths $p$ and $q$ to $Z$ in this manner. 
The notion of {\em type} applies to the vertices of the paths $p$ and $q$ so that the new 
graph $\ol{Z}$ is again bipartite. One easily sees that the bipartite graph $\ol{Z}$ 
again has girth $\ge 2n$ and that each vertex has valence $\ge 2$. 
The {\em free construction} based on a connected graph $Z_0$ of
girth $\ge 2n$ consists in the inductive application of the
bar-operation: $Z_{i+1}:= \ol{Z_i}$. Then the direct limit of the
resulting graphs is a thick building. 
We modify the above procedure by supplementing it with the operation $Z\embed Z'$ described below. 

\medskip
Let $Z$ be a bipartite graph of girth $\ge 2n$. We define 
a new graph $Z'$ as follows. For every vertex-type $l=1, 2$, every
$m\ge 3$, every $m$-tuple of mutually antipodal edges $\De_i$ 
in $Z$ and integers $0< r_i< n-1$, $i=1, ..., m$, satisfying
$$
r_i+r_j\ge n, \quad \forall i\ne j,
$$
 we attach to $Z$ an $m$-pod $T$ with the bases
$\De_1,...,\De_m$, center of the  type $l$ and the legs of the
lengths $r_1,...,r_m$ respectively. Denote the graph obtained from
$Z$ by attaching all these $m$-pods by $Z'$. Then $Z'$ is a
bipartite graph. Applying Lemma \ref{L3} repeatedly, we see that
$Z'$ still has girth $\ge 2n$.

\medskip We now proceed with the inductive construction of the building $Y$. We start
with $X_0$, which is an arbitrary connected bipartite graph of
girth $\ge 2n$.

Then set $X_1:= X_0'$ (by attaching $m$-pods for all $m$ to all $m$-tuples of pairwise antipodal chambers). 
Take $X_2:= \ol{X_1}$
(i.e, it is obtained from $X_1$ as in the free construction) and
continue this 2-step process inductively: for every even $N=2k$ 
set $X_{N+1}:= X_N'$ and $X_{N+2}:= \ol{X_{N+1}}$.

Let $Y$ denote the increasing union of the resulting graphs. Then,
clearly, $Y$ is a connected infinite bipartite graph.

\begin{lem}\label{L4}
$Y$ is a thick building modeled on $W$, satisfying Axiom A.
\end{lem}
\proof 1. Clearly, $Y$ has girth $\ge 2n$. Note that for each $N$,
the natural inclusion $X_N\to X_{N+1}$ is 1-Lipschitz 
(distance-decreasing). Moreover, by the construction, the maps $X_N\to X_{N+1}$ are $n-1$-isometric in the sense of Section \ref{metric}. 
By the construction, if $x, y\in X_N$ ($N$ is odd)
are vertices within distance $d\ge n+1$, then
$$
d_{X_N}(x, y)> d_{X_{N+1}}(x,y)
$$
(as there will be a pair of vertices within distance $n+1$ on the
geodesic $\ol{xy}\subset X_n$, while their distance in $X_{N+1}$
becomes $n-1$). Thus,
$$
d_{X_{N+2s}}(x, y)\le n, \quad \hbox{~~where~~} s=d-(n+1).
$$
Therefore, $Y$ has the diameter $n$. For every vertex $y\in Y$ there exists a vertex $y'\in Y$ which has 
the (combinatorial) distance $n$ from $y$. Therefore, attaching the $q$-paths in the bar-operation 
assures that there are infinitely many half-apartments in $Y$ connecting $y$ to $y'$. In particular, there 
are infinitely many apartments containing $y$ and $y'$. This implies that $Y$ is a thick (with each vertex having infinite valence)  
spherical building with the Weyl group $W=I_2(n)$.

\medskip 
2. In order to check Axiom A, let $\De_1,...,\De_m\subset Y$ be
antipodal chambers and $r_1,...,r_m$ be positive integers so that
$$
r_i+r_j\ge n, \quad \forall i\ne j.
$$

Then there exists $k_0$ so that $\De_1,...,\De_m\subset X_{k_0}$.
Since the maps $X_k\to Y$ are distance-decreasing, it follows that
there exists $k_1\ge k_0$ so that $\De_1,...,\De_m\subset X_{k_1}$
are antipodal. Therefore, by the construction, for every odd step
of the induction there will be two $m$-pods with the legs of the
lengths $r_1,...,r_m$ and centers of the type $l=1, 2$ attached to
the bases $\De_1,...,\De_m$. Therefore, the intersections
$$
\bigcap_{i=1}^m S_{r_i}(\De_i)\subset X_k, k\ge k_1,
$$
will contain at least $(k-k_1)/2$ vertices of both types. 
Since the maps $\iota: X_N\to X_{N+k}, k\ge 0$, are $(n-1)$-isometric, 
$\iota(S_r(\De))\subset S_r(\De)\subset X_{N+k}$. Therefore, 
$Y$ satisfies Axiom A. \qed

This concludes the proof of Theorem \ref{Aexists}. \qed

\medskip
One can modify the above construction by allowing the transfinite induction, but we will not need this. More
interestingly, one can modify the construction of $Y$ to obtain a
rank 2 spherical building $X$ which satisfies the following {\em
universality property} (with $n$ fixed): 

{\bf Axiom U.} Let $G$ be an arbitrary finite connected bipartite
graph of girth $\ge 2n$, let $H\subset G$ is a (possibly
disconnected) subgraph and $\phi: H\to X$ be a morphism (a distance-decreasing
embedding preserving the type of vertices). Then $\phi$ extends to a morphism $G\to X$.

\medskip
Thus, Axiom A is a special case of the Axiom U, defined with
respect to a particular class of graphs $G$ (i.e. $m$-pods), their
subgraphs (the sets of vertices $x_i$ of valence $1$) and maps
$\phi$ (sending $x_i$'s to vertices of antipodal chambers). The
Axiom U is somewhat reminiscent of the Kirszbraum's property (see
e.g. \cite{LPS}). With this in mind, the construction of
(countably infinite) buildings satisfying Axiom U is identical to
the proof of Theorem \ref{Aexists}.


%

\section{Highly homogeneous buildings satisfying Axiom A}

The goal of this section is to show that the
``highly homogeneous'' buildings constructed by K.~Tent in \cite{Tent},
satisfy Axiom A. This will give an alternative proof of Theorem \ref{Aexists}.

We need several definitions. From now on, fix an integer $n\ge 2$. 

\begin{defn}
Let $G$ be a finite graph $G$ with the set of vertices $V(G)$ and the set of edges $E(G)$. 
Define the {\em weighted Euler characteristic}
of $G$ as
$$
y(G)=(n-1)|V(G)| - (n-2)|E(G)|.
$$
\end{defn}

Define the class of finite graphs $\K$ as the class of bipartite graphs $G$ satisfying the following:

1) girth$(G)\ge 2n$.

2) If $G$ contains a subgraph $H$ which in turn contains an embedded $2k$-cycle, $k>2m$, then
$$
y(H)\ge 2n+2.
$$

We convert $\K$ to a category, also denoted $\K$, by declaring morphisms between graphs in $\K$ to be
label-preserving embeddings of bipartite graphs which are $n-1$-isometric maps with respect to the
combinatorial metrics on graphs.

A bipartite graph $U$ is called a $\K$-{\em homogeneous universal model} if it satisfies the following:

1) $U$ is terminal for the category $\K$, i.e.: Every finite
subgraph in $U$ belongs to $\K$ and for every graph  $G\in \K$
there exists an $(n-1)$-isometric embedding  $G\to U$.

2) If $G\in \K$ and $\phi, \psi: G\to U$ are $(n-1)$-isometric
embeddings, then there exists an automorphism
$$
\al: U\to U
$$
so that $\al\circ \phi=\psi$.

The main result of \cite{Tent} is

\begin{thm}\label{uni}
The category $\K$ admits a $\K$-homogeneous universal model $X$.
\end{thm}

Most of the proof of the above theorem deals with establishing that the category $\K$
satisfies the following {\em amalgamation (or pull-back) property}:

Every  diagram
\begin{diagram}
 &  &     X_1  &  & \\
 &\ruTo &  &  & \\
X_0 &  &  &  & \\
 &\rdTo &  &  & \\
 &  &     X_2  &  &
\end{diagram}
extends to a commutative diagram
\begin{diagram}
 &  &     X_1  &  & \\
 & \ruTo &  & \rdTo & \\
X_0 &  & \rTo  &  & X_3\\
 & \rdTo &  & \ruTo & \\
 &  &     X_2  &  &\\
\end{diagram}

The universal graph $X$ as in Theorem \ref{uni} is then shown in
\cite{Tent} to be a rank 2 thick spherical building with the Weyl
group $W=I_2(n)$, such that the automorphism group $Aut(X)$ of $X$ acts
transitively on the set of apartments in $X$, so that the
stabilizer of every apartment is infinite and contains $W$.
Moreover,  $Aut(X)$ also acts transitively on the set of simple
$2(n+1)$-cycles in $X$.

\begin{prop}
The universal graph $X$ as above satisfies Axiom A.
\end{prop}
\proof Let $T$ be an $m$-pod with the bases $\De_1,...,\De_m$ and
the legs of the length $r_1,...,r_m$. Since $T$ contains no
embedded cycles, it is an object in $\K$. In particular, the
graph $B$ which is the disjoint union of the bases
$\De_1,...,\De_m$, is also an object in $\K$. We let $\psi: B\to
X$ denote the identity embedding, which is necessarily a morphism.
Since the chambers $\De_i$ are super-antipodal in $T$ (i.e.,
$\De_i, \De_j$ are within distance $\ge n-1$ for all $i\ne j$), it
follows that the embedding $B\to T$ is a morphism in $\K$. Then,
by repeatedly using the amalgamation property, we can amalgamate
$N$ copies of $T$ along $B$ to obtain a graph $G_N$ which is again
an object in $\K$.

Let $\De_1,...,\De_m$ be a collection of antipodal chambers in $X$.
Let $r_1,...,r_m$ be numbers so that
$$
r_i+r_j\ge n, \forall i\ne j.
$$

We claim that
$$
\bigcap_{i=1}^m B_{r_i}(\De_i)
$$
contains infinitely many vertices of the given type $l=1, 2$.

Indeed, the disjoint union of the chambers $\De_i$ determines a
bipartite graph $B$. Form an $m$-pod $T$ with the union of bases
$B$, legs of the lengths $r_1,...,r_m$ and the center $z$ of the
type $l$. Let $G_N$ be the graph obtained by amalgamating $N$
copies of $T$ as above along the bases. Clearly,
$$
\bigcap_{i=1}^m B_{r_i}(\De_i)\subset G_N
$$
contains $N$ vertices of the type $l$, the centers of the
$m$-pods $T$. Since $G\in \K$, and $X$ is terminal with respect to
$\K$, it follows that  there exists an $(m-1)$-isometric embedding
$\phi: G_N\to X$. Because $\phi$ is distance-decreasing, the
intersection
$$
\bigcap_{i=1}^m B_{r_i}(\phi(\De_i))\subset X
$$
also contains $N$ vertices of the type $l$.

We thus obtain two morphisms $\phi, \psi: B\to X$, where $\psi$ is the identity embedding. By the property
3 of a $\K$-{\em homogeneous universal model}, there exists an
automorphism $\al: X\to X$ so that $\al\circ \phi=\psi$.
Therefore,
$$
\bigcap_{i=1}^m B_{r_i}(\De_i)\subset X
$$
contains at least $N$ vertices of the type $l$, namely, the images of the centers of the $m$-pods $T\subset G$ under $\al\circ \phi$. 
Since $N$ was chosen arbitrary, the proposition follows. \qed

\section{Intersections of balls in buildings satisfying Axiom
A}\label{sec:int}

In this section we prove several basic facts about cardinalities
of intersections of balls in buildings satisfying Axiom A.

\begin{lem}\label{L1}
Suppose that $X$ is a thick spherical building with the Weyl group
$W=I_2(n)$. Let $r_1+r_2=n-1$ and $\De_1, \De_2$ be non-antipodal
chambers (i.e. they are within distance $\le n-2$). Then
$B_{r_1}(\De_1)\cap B_{r_2}(\De_2)$ contains vertices of both
types.
\end{lem}
\proof Let $A\subset X$ denote an apartment containing $\De_1,
\De_2$. It suffices to consider the case when the distance between
the chambers is exactly $n-2$ (as the chambers get closer the
intersection only increases). We will assume that $r_1>0, r_2>0$ and will leave the remaining cases to the reader. 
Then $A$ will contain unique vertices $x, y$ (of distinct type) so that 
$$
d(x, \De_1)=r_1, \quad d(x, \De_2)=r_2-1, \quad d(y, \De_1)=r_1-1, \quad d(y, \De_2)=r_2. 
$$
Thus, $x, y\in B_{r_1}(\De_1)\cap B_{r_2}(\De_2)$. (Note that if $d(\De_1, \De_2)=n-2$ then $\{x, y\}= 
B_{r_1}(\De_1)\cap B_{r_2}(\De_2)$.)  \qed

\begin{lem}\label{L2}
For every  thick spherical building $X$ with the Weyl group 
$W=I_2(n)$, and every pair of antipodal chambers  $\De_1, \De_2\subset X$, 
and non-negative integers $r_1, r_2$ satisfying $r_1+r_2=n-1$, the intersection
$$
B_{r_1}(\De_1)\cap B_{r_2}(\De_2)
$$
consists of exactly two vertices, one of each type. 
\end{lem}
\proof Let $A\subset X$ be an apartment containing $\De_1, \De_2$.
It is clear that the intersection
$$
B_{r_1}(\De_1)\cap B_{r_2}(\De_2) \cap A
$$
consists of exactly two vertices $u, v$, one of each type types. 
Let $\al\subset A$ denote the subarc of length $n-1$
connecting vertices $x_i$ of the chambers $\De_i, i=1, 2$, so that $u\in \al$.  
Suppose there is a vertex $z\in X \setminus A$, $type(z)=type(u)$, so that
$$
d(z, \De_i)=r_i, i=1, 2.
$$
Then it is clear that $d(x_i, z)=r_i$, $i=1, 2$ and we thus obtain
a path (of length $n-1$)
$$
\be=\ol{x_1 z}\cup \ol{z x_2}
$$
connecting $x_1$ to $x_2$. Since $d(x_1, x_2)=n-1$, it follows that $\be$ is a geodesic path in $X$. 
Thus, we have two distinct geodesics $\al, \be\subset X$ of the length $n-1$ connecting $x_1, x_2$. The 
union $\al\circ \be $ is a (possibly nonembedded) homologically nontrivial cycle of length $2(n-1)$ in $X$. 
This contradicts the fact that $X$ has girth $2n$. \qed  


\begin{lem}\label{L0}
Let $X$ be a  building with the Weyl group
$W=I_2(n)$, satisfying Axiom A. Suppose that $r_i,
i=1,...,m$ are positive integers so that
\begin{equation}\label{e1}
r_k\le n-1, k=1,...,m,
\end{equation}
and
\begin{equation}\label{e2}
\sum_i r_i \ge (n-1)(m-1).
\end{equation}
Then for every $m$-tuple of antipodal chambers $\De_1, ... \De_m$
in $X$, one of the following mutually exclusive cases 
occurs:

a) Either the intersection
$$
\bigcap_{i} B_{r_i}(\De_i)
$$
contains infinitely many vertices of both types $l=1, 2$.

b) Or (\ref{e2}) is the equality, for two indices, $i\ne j$,
$r_i+r_j=n-1$ and for all $k\notin \{i, j\}$ the inequality
(\ref{e1}) is the equality.
\end{lem}
\proof If $r_i+r_j\ge n$ for all $i\ne j$, the assertion follows
from Axiom A (namely, the alternative (a) holds). Suppose that, say, $r_1+r_2\le n-1$. Then
$$
\sum_{i=1}^m r_i \le (n-1)+ \sum_{i=3}^m r_i\le (n-1)(m-1).
$$
Since $\sum_i r_i \ge (n-1)(m-1)$, we see that $r_1+r_2= n-1$,
$r_3=...=r_m=n-1$ and $\sum_{i=1}^m r_i =(n-1)(m-1)$. The fact that (a) and (b) cannot occur simultaneously, follows from Lemma \ref{L2}. 
\qed

\begin{cor}\label{C2}
Under the assumptions of Lemma \ref{L2}, let $\De_1,...,\De_m$ be
antipodal chambers in $X$. Then 

 1. If $r_1, ..., r_m$
are non-negative integers so that $r_1+r_2=n-1$ and
$r_3=...=r_m=n-1$, then the intersection of balls
$$
I:=\bigcap_{i} B_{r_i}(\De_i)
$$
consists of exactly two vertices (one of each type).

2. If $r_1, ... r_m$ are integers so that $\sum_i r_i <
(n-1)(m-1), 0\le r_i\le n-1, i=1, ..., m$ then the above
intersection of balls $I$ is empty.
\end{cor}
\proof The first assertion follows from Lemma \ref{L2}, since
$B_{r_i}(\De_i)=X, i\ge 3$. To prove the second assertion we note
that there are $i\ne j\in \{1, ..., m\}$ so that $r_i+r_j<n-1$.
Therefore, $B_{r_i}(\De_i)\cap B_{r_j}(\De_j)=\emptyset$ since
$d(\De_i, \De_j)=n-1$. \qed

\medskip
Recall that $X_l$ denotes the set of vertices of type $l$ in $X$. 
By combining Lemma \ref{L0} and Corollary \ref{C2}, we obtain

\begin{cor}\label{C3}
Suppose that $X$ satisfies Axiom A. Let $\De_1,...,\De_m$ be
antipodal chambers in $X$ and $r_1, ..., r_m$ are non-negative
integers so that $r_i\le n-1$, $i=1,...,m$. Then the following are
equivalent:

1) $\bigcap_{i} B_{r_i}(\De_i)\cap X_l$ is a single point for
$l=1, 2$.

2) After renumbering the indices, $r_1+r_2=n-1$ and
$r_3=...=r_m=n-1$.

Moreover, if $\sum_{i} r_i\ge (n-1)(m-1)$ then $\bigcap_{i} B_{r_i}(\De_i)$ contains vertices of both types.
\end{cor}

\section{Pre-rings}\label{under}

An {\em pre-ring} is an algebraic system $R$ with the usual
properties of a ring, except that the operations are only partially defined. 
(By analogy with groupoids, the pre-rings should be
called {\em ringoids}, however, this name is already taken for something else.)

The standard examples of pre-rings which are used in 
calculus are $\widehat{\R}=\R\cup \pm \infty$ and $\widehat{\C}=\C \cup \infty$. Below is a similar
example which we will use in this paper. For a ring $R$  define the pre-ring
$\widehat{R}:= R\cup \infty$.
The algebraic operations in $\widehat{R}$ are extended from the ring $R$
as follows:

1. Addition and multiplication are commutative and associative;
$0$ and $1$ are neutral elements with respect to the addition and
multiplication.

2. Moreover, we have
$$
\begin{array}{|c|c|c|}
\hline
\hbox{addition} & x\ne \infty &\infty\\
\hline
y\ne \infty  & x+y & \infty  \\
\hline \infty & \infty & \hbox{undefined} \\ \hline
\end{array} \quad
\begin{array}{|c|c|c|c|}
\hline
\hbox{multiplication} & 0 & x\in R\setminus \{0\} &   \infty\\
\hline
        0 &     0 &  0 & 0\\
\hline
y\in R\setminus \{0\}  & 0 & xy & \infty \\
\hline
            \infty & 0 & \infty & \infty\\ \hline
\end{array}
$$

\begin{rem}
It is customary to assume that $0\cdot \infty$ is undefined, but in the situation we are interested in (where pre-rings will appear 
as degenerations of rings), we can assume that $0\cdot \infty=0$. 
\end{rem}

\section{Schubert pre-calculus}\label{sec:calculus}

From now on, we fix a thick spherical building $X$ satisfying Axiom A. (Much
of our discussion however, uses only the fact that $X$ is a
spherical building with the Weyl group $I_2(n)$.)

Our next goal is to introduce a {\em  Schubert pre-calculus} in
$X$. According to a theorem of Kramer and Tent \cite{KT}, for $n\notin \{2, 3, 4, 6\}$, there are no
thick spherical buildings with the Weyl group $I_2(n)$ that admit 
structure of an algebraic variety defined over an algebraically closed 
field. Since we are interested in general $n\ge 2$, this forces the algebro-geometric features 
of the buildings described below to be quite limited. 

Let $l\in \{1, 2\}$ be a type of vertices of $X$. We will think of the set $X_l$ of points of type $l$ as the $l$-th ``Grassmannian''. 
Let $\De\subset X$ be a chamber and $0\le r\le n-1$ be an integer. 
We define the ``Schubert cell'' $C_{r}(\De)\subset X_l$ to be the $r$-sphere $S_r(\De)$ in $X_l$
centered at $\De$ and having radius $r$:
$$
C_{r}(\De)= \{ x\in X_l: d(x, \De)=r\}.
$$
(We suppress the dependence on $l$ in the notation for the
Schubert cell.) The number $r$ is the ``dimension'' of the cell.
We define the ``Schubert cycle'' $\ol{C_{r}(\De)}$, the
``closure'' of the Schubert cell $C_{r}(\De)$,  as the closed
$r$-ball centered at $\De$:
$$
\ol{C_{r}(\De)}:= B_r(\De)\cap X_l.
$$
The number $r$ is the ``dimension'' of this cycle. Thus, each $r$-dimensional Schubert cycle is the union of 
$r+1$ Schubert cells which are the ``concentric spheres''. By taking $r=n-1$, 
we see that $X_l$ is a Schubert cycle of dimension $n-1$. 

There is much more to be said here, but we defer this discussion to another paper.

\medskip
{\bf Homology}.  The coefficient system for our homology pre-ring is
the pre-ring $\widehat{R}$ defined in Section \ref{under}. 
The simplest case will be when $R=\Z/2$, then  $\widehat{R}$ consists of three elements: $0, 1, \infty$.
This example will actually suffice for our purposes, but our discussion here is more general.
We will suppress the coefficients in the notation for $H_*(X_l, \widehat{R})$  in what follows.

Let $W=I_2(n)$. We declare $d=n-1$ to be the formal dimension of $X_l$.
Set $r^*:= d-r$ for $0\le r\le d$. Fix a (positive) chamber $\De_+\subset X$.

Using the  Schubert pre-calculus we define the {\em homology pre-ring}
$H_*(X_l)$ ($l=1, 2$) with coefficients in $\widehat{R}$, by declaring its
(additive) generators in each dimension $0\le r\le n-1$ to be the
Schubert classes $[\ol{C_{r}(\De)}]$, where $\De$ are chambers in
$X$. We declare
$$
C_r:=[\ol{C_{r}(\De)}]=[\ol{C_{r}(\De_+)}]
$$
for every $\De$ and set
$$
H_r(X_l)=0, r<0, \quad r> d.
$$
The ``fundamental class'' in $H_{d}(X_l)$ is represented by
$X_l=B_{d}(\De_+)$. We declare a collection of cycles  $\ol{C_{r_i}(\De_i)}, i\in I$ to be {\em transversal} if the chambers
$\De_i, i\in I$ are pairwise antipodal. Using this notion of transversality we define the {\em intersection product} 
on $H_*(X_l)$ as follows.

Consider two antipodal chambers $\De_1, \De_2$. For $0\le r_1, r_2\le n-1$,
$$
\ol{C_{r_1}(\De_1)}\cap \ol{C_{r_2}(\De_2)}= B_{r_1}(\De_1)\cap
B_{r_2}(\De_2),
$$
is the ``support set'' of the product class
$$
[\ol{C_{r_1}(\De_1)}]\cdot [\ol{C_{r_2}(\De_2)}] \in H_{r_3}(X_l),
$$
where
$$
r_3^*= r_1^* + r_2^*
$$
i.e.,
$$
r_3=r_1+r_2-(n-1).
$$
The product class itself is a multiple $a\cdot
[\ol{C_{r_3}(\De_+)}]$ of the standard generator. To compute $a\in
\widehat{R}$, we declare that the classes $c=C_{r_3}$ and
$c*=C_{r_3^*}$ are ``Poincar\'e dual'' to each
other:
$$
c=PD(c^*),
$$
as their dimensions add up to the dimension $d$ of the
fundamental class. Therefore, take a chamber $\De_3$ antipodal to
both $\De_1, \De_2$: It exists by Lemma \ref{anti}.  Then $a\in \widehat{R}$ is the cardinality  of the
intersection:
$$
\ol{C_{r_1}(\De_1)}\cap \ol{C_{r_2}(\De_2)} \cap
\ol{C_{r_3}(\De_3)}.
$$

\begin{rem}
Here and in what follows we are abusing the terminology and declare cardinality of an infinite set to be $\infty$: 
This is justified, for instance, by the fact that Theorem \ref{Aexists} 
yields buildings that have countably many vertices and 
our convention amounts to $\aleph_0=\infty\in \widehat{R}$.  
\end{rem}

As we will see below, this cardinality is $0, 1$ or $\infty$, these cardinalities are 
naturally identified with the elements of $\widehat{R}$.

One can easily check (see below) that $a$ does not depend on the
choice of cycles representing the given homology classes. In
particular, the fundamental class is the unit in the pre-ring
$H_*(X_l)$.

We now compute $a$ using the results of Section \ref{sec:A}:

1. If $r_1+r_2>n-1$ then $a=0$. (Corollary \ref{C2}, Part 2.)

2. If $r_1+r_2=n-1$ then $a=1$: The Schubert cycles
$\ol{C_{r_i}(\De_i)}, i=1, 2$, are Poincar\'e dual to each other.
(Lemma \ref{L2}.)

3. Suppose now that $r_1+r_2<n-1$.

3a. If $0<r_1, r_2<n-1$ then $a=\infty$ unless $r_1+r_2=n-1$.
(This immediately follows from Lemma \ref{L0}.) Thus, only
Poincar\'e dual classes have ``finite'' intersection.

3b. If $r_i=n-1$ or $r_i=0$ for some $i=1, 2$, then $a=1$.
(Obvious.)

\begin{lem}
Let $C_{r_i}\in H_{r_i}(X_l), i=1,...,m$ be the generators (the
Schubert classes) so that
$$
r_1+...+r_n=d(m-1)  \iff \sum_{i=1}^m r_i^*=d,
$$
i.e., the product of these classes (in some order) equals $a[pt]$,
where $pt=\ol{C_0(\De_+)}$. Then $a\in \widehat{R}$ is the cardinality of
the intersection
$$
\bigcap_{i=1}^m B_{r_i}(\De_i),
$$
where $\De_1,...,\De_m$ are pairwise antipodal chambers (which exist by Lemma \ref{anti}).
\end{lem}
\proof First of all, without loss of generality we may assume that
none of the classes $C_{r_i}$ is the unit $[X_l]$ in $H_*(X_l)$. Note
that, since $r_1+...+r_m=(n-1)(m-1)$, in the computation of the
product of $C_{r_1},...,C_{r_m}$ we will never encounter the
multiplication by zero. Then (after permuting the indices), the
product of the classes $C_{r_1},...,C_{r_m}$ will be of the form
$$
... (C_{r_1} \cdot C_{r_2}) ...
$$
By the definition, $C_{r_1} \cdot C_{r_2}=a_{12} C_r$, where
$$
r^*=r_1^*+r_2^*.
$$
The element $a_{12}\in \widehat{R}$ is the cardinality of the intersection
$$
B_{r_1}(\De_1)\cap B_{r_2}(\De_2)\cap B_{r^*}(\De),
$$
where $\De_1, \De_2, \De$ are pairwise antipodal. In view of the
above product calculations 1---3, and the fact that $r_1\ne n-1,
r_2\ne n-1$, we see that $a_{12}=\infty$ (since $a_{12}=0$ is
excluded), unless $r=d$, $r_1=r_2^*$ and, therefore,
$c_2=PD(c_1)$. In the latter case, $a_{12}$ is the cardinality
(equal to $1$) of the intersection
$$
B_{r_1}(\De_1)\cap B_{r_2}(\De_2).
$$
Since
$$
\sum r_i=(n-1)(m-1), \quad 0\le r_i\le n-1, \quad i=1,...,m,
$$
we conclude that $r_3=...=r_m=n-1$. Thus, $m=2$ and $a=a_{12}=1$
in this case.

If $a_{12}=\infty$ then it follows from the definition of $\widehat{R}$
that $a=\infty$, since, in the computation of the product of
 $C_{r_1},...,C_{r_m}$ we will never multiply by zero. On the other hand,
in this case the classes $C_{r_1}, C_{r_2}$ are not Poincar\'e dual to
each other and Lemma \ref{L0} implies that the intersection
$$
\bigcap_{i=1}^m B_{r_i}(\De_i)\subset X_l
$$
is also infinite. Lemma follows. \qed

\begin{cor}
$H_*(X_l, \widehat{R})$ is a pre-ring.
\end{cor}
\proof The only thing which is unclear from the definition is that
the product is associative. To verify associativity, we have to
show that
\begin{equation}\label{asso}
((C_{r_1} C_{r_2}) C_{r_3})\cdot C_{r_4}= (C_{r_1} (C_{r_2}  C_{r_3}))\cdot C_{r_4}
\end{equation}
where $C_{r_i}\in H_{r_i}(X_l)$ are the generators and
$$
r_1+r_2+ r_3+ r_4= (4-1)(n-1).
$$
However, the equality (\ref{asso}) immediately follows from the
above lemma. \qed

\medskip 

Similarly to the definition of the Schubert pre-calculus on the Grassmannians $X_l$, we define the Schubert pre-calculus on 
the ``flag-manifold'' $Fl(X)$ associated with $X$, i.e., the set of edges $E(X)$ of the graph $X$ underlying the building $X$. 
The set $E(X)$ will be identified with the set of mid-points of the edges. We have two projections
$$
p_l: E(X)\to X_l, l=1, 2
$$
sending each edge to its end-points. We will think of these projections as ``$\P^1$-bundles.'' Accordingly, we define Schubert cycles in $Fl(X)$ 
by pull-back of Schubert cycles in $X_l$ via $p_l$:
$$
\ol{C_{r,l}(\De)}:=p_l^{-1}\left(\ol{C_{r-1}(\De)}\right), r=1,...,n. 
$$
while $0$-dimensional cycles in $Fl(X)$ are, of course, just the edges of $X$. In terms of metric geometry of $X$, 
the cycles  $\ol{C_{r+1,l}}(\De)$ are described as follows. Fix a chamber $\De$. Define the Schubert cell $C_{r,l}(\De)$ 
to be the set of chambers $\De'\subset X$ so that the distance between the midpoints $mid(\De), mid(\De')$ 
of $\De, \De'$ equals $r$ and the minimal distance $r-1$ between $\De, \De'$ is realized 
by a vertex of type $l$ in $\De'$. Here the convention is that $C_{r,l}(\De)=C_{r,l+1}(\De)$ for $r=0, r=n=girth(X)/2$, since for these 
values of $r$ the minimal distance is realized by vertices of both types. The corresponding Schubert cycles $\ol{C_{r,l}(\De)}$ are 
defined by adding to  $C_{r,l}(\De)$ all the chambers $\De''$ contained in the geodesics connecting 
$mid(\De), mid(\De')$, for $\De'\in C_{r,l}(\De)$. The notions of transversality as in the case of $X_l$, is given by taking antipodal chambers. 
The Poincar\'e Duality is defined by 
$$
PD([\ol{C_{r,l}(\De)}])=[\ol{C_{n-r,3-l}(\De)}], l=1, 2. 
$$
The reader will verify that this is consistent with the property that the intersection
$$
\ol{C_{r_1,l}(\De_1)}\cap \ol{C_{r_2,3-l}(\De_2)}
$$
is a single point. We declare that the homology classes  $[\ol{C_{r,l}(\De)}]$ are independent of $\De$ and 
set up the notation 
$$C_{r,l}:=C_{w}:=[\ol{C_{r,l}(\De)}],$$
 where $w\in W$ is the unique element such that 
$w(\De)\in C_{r,l}(\De)$.  Then the Poincar\'e Duality takes the form
$$
PD(C_w)=C_{w_\circ w}, 
$$
where $w_\circ\in W$ is the longest element. 

We  declare that $C_{r,l}$, $r=0,...,n, l=1, 2$, form a basis of  $H_*(Fl(X))$, where $r=\dim(C_{r,l})$. 
We also require the pull-back maps $p_l$ to be pre-ring homomorphisms. 
It remains to define the intersection products of the form
$$
C_{r_1,1}\cdot C_{r_2,2}, \quad 0\le r_1, r_2\le n. 
$$
Analogously to the product in $H_*(X_l)$, we take two antipodal chambers $\De_1, \De_2$ and set
$$
C_{r_1,1}\cdot C_{r_2,2}= a_1 C_{r_3, 1}+ a_2 C_{r_3, 2}, \quad a_l\in \widehat{R}, l=1, 2. 
$$
In order to compute $a_l$'s we take the third chamber $\De_3$ antipodal to $\De_1, \De_2$, $r_3:=r_1+r_2-n$, and let $a_l$ denote the cardinality of 
the intersection 
$$
\ol{C_{r_1,1}(\De_1)} \cap \ol{C_{r_2,2}(\De_2)} \cap  \ol{C_{r_3,3-l}(\De_3)}. 
$$
With these definition, we obtain a homology pre-ring $H_*(Fl(X), \widehat{R})$ abbreviated to $H_*(X, \widehat{R})$ or even $H_*(X)$. 
The proof of the following proposition is similar to the case of $H_*(X_l)$ and is left to the reader:

\begin{prop}\label{H_*(X)}
Let $X$ be a thick building with the Weyl group $I_2(n)$, satisfying Axiom A. Then $H_*(X)$ is an associative and commutative 
pre-ring, generated by the elements $c_{l,r}, l=1, 2, r=0,...,n$, subject to the relations:
\begin{enumerate}
\item $$
C_{1,0}=C_{2,0}, 
$$
$$
C_{1,n}=C_{2,n}=1,
$$
is the unit in $H_*(X)$, 
\item 
$$
C_{r_1,l}\cdot C_{r_2,l}=0, \hbox{~~if~~} r_1+r_2\le n, \quad l=1, 2, 
$$ 
\item 
$$
C_{r_1,l}\cdot C_{r_2,l}=\infty, \hbox{~~if~~} n<r_1+r_2, \quad l=1, 2,  
$$
\item 
$$
C_{r_1,1}\cdot C_{r_2,2}= 1, \hbox{~~if~~} r_1+r_2=n,
$$
\item 
$$C_{r_1,1}\cdot C_{r_2,2}=0, \hbox{~~if~~} r_1+r_2<n,$$ 
\item 
$$
C_{r_1,1}\cdot C_{r_2,2}= \infty C_{r_3, 1}+ \infty C_{r_3, 2}, \hbox{~~ if~~} r_1+r_2>n, \hbox{~~  where~~} r_3=(r_1+r_2)-n. 
$$
\end{enumerate}
\end{prop}

\section{The stability inequalities}
\label{sec:sineq}

Suppose that $X$ is a rank 2 thick spherical building with the
Weyl group $W\cong I_2(n)$, satisfying Axiom A. We continue with the
notation from Section \ref{sec:calculus}. Recall that $\angle$ is
a path metric on $X$ so that the length of each chamber is
$\pi/n$.

We start with few simple observations. Let $C_r(\De)$ be a
Schubert cell in $X_l$ and $\eta\in C_r(\De)$, i.e., $d(\eta,
\De)=r$. Then the point $\zeta$ in $\De$ nearest to $\eta$ has the
type $l+ r$ (mod 2). In particular, $\zeta$ depends only on the
cell $C_r(\De)$ (and not on the choice of $\eta$ in the cell). Let
$\xi\in \De$ be a point within $\angle$-distance $\tau$ from
$\zeta$. Then
$$
\angle(\eta, \xi)= r\frac{\pi}{n} + \tau.
$$
In particular, this angle is completely determined by the angle
$\tau$, by the type of $\eta$ and the fact that we are dealing
with the Schubert cell $C_r(\De)$. In particular, it follows that
for each $\eta\in \ol{C_{r-1}(\De)}=\ol{C_r(\De)}\setminus
C_r(\De)$, we have
$$
\angle(\eta, \xi)< r\frac{\pi}{n} + \tau,
$$
where $\xi$ is defined as above. We now introduce the following
system of inequalities $WTI$ (weak triangle inequalities) on
$m$-tuples of vectors $\ov\la=(\la_1,...,\la_m)=(\mu_1\xi_1,...,
\mu_m\xi_m)$, with $\xi_i\in \De_+$ and $\mu_i\in \R_+$.

Each Grassmannian $X_l$ (or, equivalently, the choice of a vertex
$\zeta$ of the standard spherical chamber $\De_+$)  will contribute a
subsystem $WTI_l$ of the triangle inequalities. Consider all
possible $m$-tuples $(w_1,...,w_m)$ of elements of $W$, so
that all but two $w_i$'s are equal to $w_\circ$ (the longest element
of $W$)  and the remaining elements $w_i, w_j$ are ``Poincar\'e
dual'' to each other ($w_i=PD_l(w_j)$), i.e., their relative lengths
$r_i=\ell_l(w_i), r_j=\ell_l(w_j)$ in $W/W_l$ satisfy
$$
r_i+ r_j=n-1.
$$
In other words, the corresponding Schubert cycles
$$
C_{r_i}=[\ol{C_{r_i}(\De_+)}], \quad C_{r_j}=[\ol{C_{r_j}(\De_+)}]
$$
in $X_l$ have complementary dimensions and thus are Poincar\'e
dual to each other:
$$
C_{r_i}= PD( C_{r_j} ).
$$
See Section \ref{sec:calculus}. Equivalently, we are considering
$m$-tuples of integers $0\le r_k\le m-1$, which, after permutation
of indices, have the form
$$
(r_1,...,r_m)=(r_1, r_2=n-1-r_1, n-1,...,n-1).
$$
Note that $\ell_l(w_\circ)=n-1$, thus $\ell_l(w_i)=r_i, i=1,...,m$.

Lastly, for every such tuple $\ov{w}=(w_1,...,w_m)=(w_\circ,...,w_\circ,
w_i,..., PD(w_i),..., w_\circ)$ we impose on the vector $\ov{\la}$ the
inequality
\begin{equation}\label{ineq}
\sum_j \< \la_j, w_j(\zeta)\> = \sum_j \mu_j\cdot\cos\angle(\la_j, w_j(\zeta)) \le 0
\end{equation}
denoted $WTI_{l,\ov{w}}$. The collection of all these inequalities
constitutes the system of inequalities $WTI$.

\begin{thm}\label{main} 
For any rank 2 thick spherical building $X$ satisfying Axiom A with the
Weyl group $I_2(n)$, one has:

(i) The cone $\K_m(X)$ is cut out by the inequalities $WTI$.  

(ii) Moreover, if $\ov\la\in \K_m(X)$, then
there exists a semistable weighted configuration $\psi$ of the type
$\ov\la$ so that the points $\xi_i, i=1,...,m$, belong to mutually antipodal chambers.
\end{thm}
\proof Our proof essentially repeats the one in \cite[Theorem
3.33]{KLM1}. We present it here for the sake of completeness.

1 (Existence of a semistable configuration). We begin by taking a
collection of chambers $\De_1,...,\De_m\subset X$ in ``general
position,'' i.e., they are mutually antipodal. (In \cite{KLM1} one
instead takes a generic configuration of Schubert cycles in the
generalized Grassmannian, representing the given homology
classes.) Then for each $i=1,...,m$ we place the weight $\mu_i$ at
the point $\xi_i'\in \De_i$ that has the same type as $\xi_i$. We
claim that the resulting weighted configuration $\psi$ in $X$ is
semistable. Suppose not. Then, according to ``Harder-Narasimhan
Lemma'' \cite[Theorem 3.22]{KLM1}, there exists $l\in \{1, 2\}$ so that
in the Grassmannian $X_l$ there exists a unique point $\eta$ with 
the minimal (negative) slope with respect to $\psi$:
$$
slope_\psi(\eta)=-\sum_i \mu_i \cos(\angle(\eta, \xi_i'))<0,
$$
i.e.,
$$
\sum_i \mu_i \cos(\angle(\eta, \xi_i')) >0.
$$
Consider the Schubert cells
$$
C_{r_i}(\De_i), \quad i=1,...,m,
$$
where $r_i=d(\De_i, \eta)$ is the (combinatorial) distance between
the chamber $\De_i$ and the vertex $\eta\in X_l$. Thus,
$$
\eta\in J=\bigcap_{i=1}^m C_{r_i}(\De_i) \subset
\bar{J}=\bigcap_{i=1}^m B_{r_i}(\De_i)\subset X_l.
$$

By the observations in the beginning of this section, the function
$slope_\psi$ is constant on $J$. Since $slope_\psi$ attains unique
minimum on $X_l$, it follows that $J=\{\eta\}$. Moreover, if
$$
\eta'\in \bar{J}\setminus J,
$$
then
$$
slope_\psi(\eta')=-\sum_i \mu_i \cos(\angle(\eta', \xi_i'))<
slope_\psi(\eta),
$$
which contradicts minimality of $\eta$. Therefore, the
intersection  $\bar{J}$ is the single point $\eta$. Thus, the
product in $H_*(X_l)$ of the Schubert classes
$[\ol{C_{r_i}(\De_i)}]$, $i=1,...,m$, is $[pt]$ and the latter
occurs exactly when (after permuting the indices) the $n$-tuple
$(r_1,...,r_m)$ has the form
$$
(r_1,...,r_m)=(r_1, r_2=r_1^*, n-1,...,n-1),
$$
see Corollary \ref{C3}. Let $\ov{w}=(w_1, w_2, w_3,...,w_m)=(w_1,
w_\circ w_1, w_\circ,...,w_\circ)$ be the corresponding tuple of elements of
the Weyl group $W$. Note that
$$
\angle(\xi'_k, \eta)= \angle(\xi_k, w_k(\zeta))
$$
since $\eta\in C_{r_k}(\De_k)$ and $w_k(\zeta)\in C_{r_k}(\De_+)$, $k=1,...,m$. Therefore,
$$
0> slope_\psi(\eta)=-\sum_i \mu_i \cos(\angle(\eta, \xi_i'))= -\sum_i \mu_i \cos(\angle(\xi_k, w_k(\zeta))).
$$
The inequality $WTI_{l,\ov{w}}$ however requires that
$$
\sum_i \mu_i \cos(\angle(\xi_i, w_i(\zeta)))\ge 0.
$$
Contradiction. Therefore, $\psi$ is a semistable configuration.

\medskip
2. Suppose that $\psi=(\mu_1\xi'_1,...,\mu_m\xi'_m)$ is a weighted
semistable configuration in $X$ of the type
$$
\ov\la=(\mu_1\xi_1,...,\mu_m\xi_m).
$$
Consider an $m$-tuple $\ov{w}=(w_1, w_2, w_3,...,w_m)=(w_1,
w_\circ w_1, w_\circ,...,w_\circ)$ of elements of $W$ as in the definition of the inequalities $WTI$ 
(after permuting the indices we can assume that the tuple has this form). We will show 
that $\ov\la$ satisfies the inequality  $WTI_{l,\ov{w}}$ for $l=0,
1$. Fix $l$ and let $\zeta\in \De_+$ denote the vertex of type
$l$. Let $r_1,...,r_m$ be the relative lengths of $w_1,...,w_m$ in
$W/W_{l}$. Let $\De_i\subset X$ denote a chamber containing
$\xi_i'$. Note that $\ol{C_{r_i}(\De_i)}=X_l$ for each $i\ge 3$
since $r_i=n-1$. According to Lemmata \ref{L1}, \ref{L2}, the
intersection
$$
\bigcap_{k=1}^m \ol{C_{r_k}(\De_k)}= \ol{C_{r_1}(\De_1)} \cap \ol{C_{r_2}(\De_2)} \subset X_l
$$
contains a vertex $\eta\in X_l$ (possibly non-unique since $\De_1,
\De_2$, a priori, need not be antipodal). Therefore,
$$
d(\eta,\De_i)\le r_i= d(w_i(\zeta), \De_+), \quad i=1,...,m.
$$
Accordingly,
$$
\angle(\eta, \xi_i')\le \angle(w_i(\zeta), \xi_i) , \quad
i=1,...,m
$$
since $\xi_i=\theta(\xi_i')\in \De_+$. Therefore,
$$
0\le slope_\psi(\eta)= -\sum_i \mu_i\cos(\angle(\eta, \xi_i')) \le -
\sum_i \mu_i \cos(\angle(w_i(\zeta), \xi_i)),
$$
and
$$
\sum_i \mu_i \cos(\angle(w_i(\zeta), \xi_i))\le 0
$$
and, thus, $\ov\la$ satisfies  $WTI_{l,\ov{w}}$.  \qed

\begin{cor}
Theorem \ref{main}(i) holds for {\bf all} 1-dimensional 
thick spherical buildings (not necessarily satisfying Axiom A) with the Weyl group $I_2(n)$. 
\end{cor}
\proof We consider two thick spherical buildings $X, X'$, where
$X$ satisfies Axiom A. According to   Theorem \ref{T1}, $\K_m(X)=\K_m(X')$. Corollary follows from 
Theorem \ref{main} and existence of buildings satisfying Axiom A. \qed

\medskip
We now convert the system of weak triangle inequalities $WTI$ to
the form which appears in Theorem \ref{main0}. For
$$
\ov{w}=(w_1,...,w_m)=(w_1,w_\circ w_1, w_\circ,...,w_\circ),
$$
and $\la_i=m_i\xi_i, i=1,...,n$, we set $w:=w_1^{-1}$. Then, for $i\ge 3$,
$$
\< \la_i, w_i(\zeta)\>= \< w_i^{-1} \la_i, \zeta\>= \< w_\circ \la_i, \zeta\>= -\<\la_i^*, \zeta\>,
$$
while
$$
\< \la_2, w_2(\zeta)\>= \< w_2^{-1}(\la_2), \zeta\>= -\< w(\la_2^*), \zeta\>,
$$
$$
\< \la_1, w_1(\zeta)\>=\<w(\la_1), \zeta\>.
$$
Therefore, the inequality
$$
\sum_j \< \la_j, w_j(\zeta)\> \le 0
$$
is equivalent to
$$
\<w(\la_1), \zeta\> - \<w(\la_2^*), \zeta\>\le \<\sum_{j=3}^m \la^*_j,
\zeta\>.
$$
Since these inequalities hold for both vertices $\zeta$ of
$\De_+$, we obtain
$$
w(\la_1-\la_2^*)\le_{\De^*} \sum_{j=3}^m \la^*_j,\quad w\in W.
$$
This proves Theorem \ref{main0}. \qed

\begin{cor}
Let $X$ be a thick spherical building. Then the stability cone $\K_m(X)$ is a convex polyhedral cone.
\end{cor}
\proof It suffices to consider the case when $X$ does not have a spherical factor, i.e., its Coxeter complex $(S,W)$ is essential: 
$W$ has no global fixed points in $S$. 
The assertion of the corollary was proven in \cite{KLM1, KLM2} for all thick spherical buildings $X$ 
with the {\em crystallographic} Weyl group $W$, i.e., $W$ appearing as Weyl groups of complex semisimple Lie groups.  
If $W=W_1\times ... \times W_k$  is a finite Coxeter group (with $W_i$ Coxeter groups with connected Dynkin diagrams) which is a Weyl group 
of a thick spherical building $X$, then each $W_i$ is either crystallographic or is a finite dihedral group $I_2(n)$,  see \cite{Tits}. 
It is immediate from the definition of semistability that
$$
\K_m(X)=\K_m(X_1)\times ... \times \K_m(X_k),
$$
where $X_1,...,X_k$ are {\em irreducible factors of $X$} with respect to its joint decomposition into irreducible spherical subbuildings: 
The $X_i$'s are thick irreducible spherical buildings with essential Coxeter complexes and Weyl groups $W_i$, $i=1,...,k$. 
It therefore follows from the above result of \cite{KLM1, KLM2} and Theorem \ref{main} that each $\K_m(X_i)$, and, hence, $\K_m(X)$, 
is a convex polyhedral cone. \qed

\section{The universal dihedral cohomology algebra $A_t$}
\label{sect:dihedral cohomology}

In this section we construct a family of algebras $A_t$, $t\in \C^\times$ as a universal deformation of cohomology ring of 
the flag variety for each rank $2$ complex Kac-Moody group $G$ (including Lie groups $G=SL_3, Sp_4, G_2$).
 It turns out that the complexification $\C\otimes A_t$ is isomorphic to the coinvariant algebra of the dihedral group 
 $W=W_t$ acting on $\C^2$ with the parameter $t$, i.e., $t^2+ t^{-2}$ is the trace of the generator of the maximal normal cyclic subgroup of $W$.

For each integer $k\ge 0$ define the  {\it $t$-integer}  $[k]_t$ by
$$[k]_t:= \frac{t^k- t^{-k}}{t- t^{-1}}=t^{1-k}+t^{3-k}+\cdots +t^{k-3}+t^{k-1} \ .$$
It is well-known (and easy to see) that for $k,\ell\ge 0$ one has
\begin{equation}
\label{eq:multiplication of t-integers}
[k]_t[\ell]_t=[|k-\ell|+1]_t+[|k-\ell|+3]_t+\cdots+[k+\ell-1]_t \ .
\end{equation}

Now define the $t$-factorials $[m]_t!:=[1]_t[2]_t\cdots [\ell]_t$ and is the $t$-binomial coefficients by:
$$\displaystyle{{m\brack k}_t=\frac{[m]_t!}{[k]_t![m-k]_t!}} \ .$$

Note that, as the usual binomials, $t$-binomials ${m\brack k}_t$ extend naturally to $k\in \N$ and $m\in \R_+$, 
although we will  use them only for $m, k\in \Z$ unless we state otherwise.  
The $t$-binomial coefficients satisfy the symmetry 
$$
\displaystyle{{n\brack k}_t={n\brack n-k}_t}$$
and the Pascal recursion:
$$\displaystyle{{m\brack k}_t=t^k{m-1\brack k}_t+t^{k-m}{m-1\brack k-1}_t}$$

\begin{proposition}
\label{pr:binomial integer}
Each $t$-binomial coefficient ${n\brack k}_t$ belongs to $\Z[t+t^{-1}]$.

\end{proposition}

\begin{proof} We need the following result.

\begin{lemma} For all $k,\ell\ge 0$ we have
\begin{equation}
\label{eq:symmetric powers}
{\ell+k\brack k}_t=\sum_{0\le m\le k\ell} c_m \cdot [m+1]_t
\end{equation}
where each $c_m\in \Z_{\ge 0}$.
\end{lemma}

\begin{proof} Let $V_1=\C^2$ be the natural $SL_2(\C)$-module.
Denote $V_\ell=S^\ell V_1$ so that $\dim V_\ell=\ell+1$. Clearly, each $V_\ell$ is a simple module. 
For each $k\ge 0$ let  $V_{\ell,k}=S^k V_\ell$  and $\dim V_{\ell,k}=\binom{\ell+k}{k}$.
Recall that for each finite-dimensional $SL_2(\C)$-module $V$ the character $ch(V)$ is a function of $t\in \C^\times$
defined by
$$ch(V)=Tr\left(
\begin{pmatrix}
t & 0 \\
0 & t^{-1} \\
\end{pmatrix}| V\right ) \ .$$
It is easy to see that $ch(V_\ell)=[\ell+1]_t$ and $ch(V_{\ell,k})={\ell+k\brack k}_t$.
Using the decomposition of $V_{\ell,k}$ into simple $SL_2$-modules:
$$V_{\ell,k}= \sum_{0\le m\le k\ell} c_m \cdot V_m$$
where each $c_m\in \Z_{\ge 0}$ and applying $ch(\cdot)$ to it, we obtain \eqref{eq:symmetric powers}.
Lemma follows. 
\qed  \end{proof}

Furthermore, the obvious recursion $[m+1]_t=[2]_t[m]_t-[m-1]_t$,
which is a  particular case of \eqref{eq:multiplication of t-integers}, proves (by induction) that each 
$t$-number $[m+1]_t$ belongs to $\Z[t+t^{-1}]=\Z[[2]_t]$.

Combining this observation with \eqref{eq:symmetric powers}, we finish the proof of the proposition. \qed 
\end{proof}

Let $A'$ be the algebra over $\C(t)$  generated by $\sigma_1,\sigma_2$ subject to the relations
$$\sigma_1\sigma_2=\sigma_2\sigma_1,~(\sigma_1-t\sigma_2)(\sigma_1-t^{-1}\sigma_2)=0\ .$$

It is convenient to rewrite the second relation as:
\begin{equation}
\label{eq:quadratic coinvariant}
 [2]_t\sigma_1\sigma_2=\sigma_1^2+\sigma_2^2 \ .
\end{equation}

\begin{lemma}
\label{le:k-th power}
The following relations hold in $A'$:
\begin{equation}
\label{eq:k-th power}
[k+\ell]_t \sigma_1^{k}\sigma_2^{\ell}=[k]_t\sigma_1^{k+\ell}+[\ell]_t\sigma_2^{k+\ell}
\end{equation}
for all $k,\ell\ge 0$. In particular, the monomials $\sigma_i^k$, $i\in \{1,2\}$, $k\ge 0$ form a $\C(t)$-linear basis of $A'$.

\end{lemma}
\begin{proof} We proceed by induction in $\min(k,\ell)$. Indeed, if $k=0$ or $\ell=0$, we have nothing to prove. Otherwise, using \eqref{eq:quadratic coinvariant} and the inductive hypothesis, we obtain:
$$[k+\ell]_t \sigma_1^{k}\sigma_2^{\ell}=[k+\ell]_t (\sigma_1\sigma_2)\sigma_1^{k-1}\sigma_2^{\ell-1}=\frac{[k+\ell]_t}{[2]_t} (\sigma_1^2+\sigma_2^2)\sigma_1^{k-1}\sigma_2^{\ell-1}$$
$$=\frac{[k+\ell]_t}{[2]_t} (\sigma_1^{k+1}\sigma_2^{\ell-1}+\sigma_1^{k-1}\sigma_2^{\ell+1})$$
$$=\frac{1}{[2]_t} ([k+1]_t\sigma_1^{k+\ell}+[\ell-1]_t\sigma_2^{k+\ell}+[k-1]_t\sigma_1^{k+\ell}+[\ell+1]_t\sigma_2^{k+\ell})$$
$$=\frac{1}{[2]_t} ([k+1]_t+[k-1]_t)\sigma_1^{k+\ell}+\frac{1}{[2]_t}([\ell-1]_t+[\ell+1]_t)\sigma_2^{k+\ell}$$
$$=[k]_t\sigma_1^{k+\ell}+[\ell]_t\sigma_2^{k+\ell}$$
by \eqref{eq:multiplication of t-integers}.

Furthermore, the relations \eqref{eq:k-th power}  guarantee that the monomials  $\sigma_i^k$, $i\in \{1,2\}$, $k\ge 0$ span $A'$. To verify their linear independence,  let us compute the Hilbert series $h(A',z)$ of $A'$.
Clearly, the Hilbert series of the polynomial algebra $\C(t)[\sigma_1,\sigma_2]$ is  $\frac{1}{(1-z)^2}$ and the Hilbert series of any principal ideal $I$ in $\C(t)[\sigma_1,\sigma_2]$ generated by a quadratic polynomial is $\frac{z^2}{(1-z)^2}$. Therefore, the Hilbert series of the quotient algebra $\C(t)[\sigma_1,\sigma_2]/I$ is
$$\frac{1}{(1-z)^2}-\frac{z^2}{(1-z)^2}=\frac{1+z}{1-z}=1+\sum_{k\ge 1} 2z^k \ .$$
Applying this to our algebra $A'=\C(t)[\sigma_1,\sigma_2]/\langle(\sigma_1-t\sigma_2)(\sigma_1-t^{-1}\sigma_2)\rangle$ we see that each 
graded component of $A'$ is $2$-dimensional, which verifies the linear independence of the monomials. The lemma is proved. \qed 
\end{proof}

Denote $\displaystyle{\sigma^{[k]}_i:=\frac{1}{[k]_t!}\sigma_i^k}$, $i=1,2$, $k\ge 0$ the divided powers of $\sigma_i$, $i=1,2$.  Denote by $A$ the subalgebra of $A'$ 
generated over $\Z[t+t^{-1}]$ by all $\sigma^{[k]}_i$, $i\in \{1,2\}$, $k\ge 0$.

\begin{proposition}
\label{pr:k-th divided power1}
The following relations hold in $A$:
\begin{equation}
\label{eq:k-th divided power1}
\sigma_1^{[k]}\sigma_1^{[\ell]}={k+\ell\brack k}_t\sigma_1^{[k+\ell]}, ~\sigma_2^{[k]}\sigma_2^{[\ell]}={k+\ell\brack k}_t\sigma_2^{[k+\ell]} \ ,
\end{equation}
for all $k,\ell\ge 0$.
\begin{equation}
\label{eq:k-th divided power2}
\sigma_1^{[k]}\sigma_2^{[\ell]}={k+\ell-1\brack k-1}_t\sigma_1^{[k+\ell]}+{k+\ell-1 \brack \ell-1}_t\sigma_2^{[k+\ell]}
\end{equation}
for all $k,\ell\ge 0$.

In particular, monomials $\sigma_i^{[k]}$, $i=1,2$, $k\ge 0$, form a $\Z[t+t^{-1}]$-linear basis in $A$, and the relations \eqref{eq:k-th divided power1} 
and \eqref{eq:k-th divided power2} are defining for $A$.
\end{proposition}

\begin{proof} We have $\sigma_i^{[k]}\sigma_i^{[\ell]}=\frac{1}{[k]_t![\ell]_t!}\sigma_i^{k+\ell}=\frac{[k+\ell]_t}{[k]_t![\ell]_t!}\sigma_i^{[k]}\sigma_i^{[\ell]}$ 
for $i\in \{1,2\}$, $k\ge 0$,
which verifies \eqref{eq:k-th divided power1}. Furthermore, \eqref{eq:k-th power} implies that
$$\sigma_1^{[k]}\sigma_2^{[\ell]}=\frac{1}{[k]_t![\ell]_t!}\sigma_1^k\sigma_2^\ell=\frac{1}{[k]_t![\ell]_t![k+\ell]_t}([k]_t\sigma_1^{k+\ell}+[\ell]_t\sigma_2^{k+\ell})$$
$$=\frac{[k+\ell-1]_t!}{[k-]_t![\ell]_t!}\sigma_1^{[k+\ell]}+\frac{[k+\ell-1]_t!}{[k]_t![\ell-1]_t!}\sigma_2^{[k+\ell]} \ ,$$
which verifies \eqref{eq:k-th divided power2}.

Since all structure constants of $A$ are $t$-binomial coefficients,  Proposition \ref{pr:binomial integer} 
guarantees that $A$ is defined over $\Z[t+t^{-1}]$.

Since, as a $\Z[t+t^{-1}]$-module, $A$ is spanned  by all products of various $\sigma_i^{[k]}$ and each such a monomial is a 
$\Z[t+t^{-1}]$-linear combination of divided powers $\sigma_i^{[\ell]}$, $i\in \{1,2\}$, $\ell\ge 0$ by
\eqref{eq:k-th divided power1} and \eqref{eq:k-th divided power2}, we see that the divided powers span $A$  
as a $\Z[t+t^{-1}]$-module. It is also clear that  the divided powers $\sigma_i^{[\ell]}$, $i\in \{1,2\}$, $\ell\ge 0$ are 
$\Z[t+t^{-1}]$-linearly independent because that was the case in $A'$ by Lemma \ref{le:k-th power}. 
Therefore, relations \eqref{eq:k-th divided power1} and \eqref{eq:k-th divided power2} are defining.
The proposition is proved. \qed \end{proof}

Now we will use the standard algebraic trick of specializing a formal parameter $t$ into a non-zero complex number $t_0$. 
Clearly, this is impossible to do for $A'$ because it is defined over $\C(t)$ but is a perfectly reasonable to do so for the 
algebra $A$ which is defined over $\Z[t+t^{-1}]$. Indeed for each $t_0\in \C^\times$ we define $\tilde A_{t_0}=R_0\otimes_R A$, 
where $R=\Z[t+t^{-1}]$, $R_0=\Z[t_0+t_0^{-1}]\subset \C$, where $R_0$ is regarded as an $R$-module via the evaluation homomorphism 
$R\to R_0$ which takes $t$ to $t_0$. By the construction, $\tilde A_{t_0}$ is a free $\Z[t_0+t_0^{-1}]$-module, e.g., 
it has a basis $\sigma_i^{[k]}$, $i\in\{1, 2\}$, $k\ge 0$.

With a slight abuse of notation, from now on we will denote by $t$ a non-zero complex number so that $\tilde A_t$, $t\in \C^\times$  
is the family of unital $\Z[t+t^{-1}]$-algebras with the
presentation \eqref{eq:k-th divided power1} and \eqref{eq:k-th divided power2} (and $\sigma_1^{[0]}=\sigma_2^{[0]}=1$).

For each $t\in \C^\times\setminus\{-1,1\}$ define $n_t\in \Z\sqcup \{\infty\}$ to be the order of $t^2$ in the multiplicative group 
$\C^\times$. If $t=\pm 1$, we set $n_{\pm 1}:=\infty$. 
Thus, $n_t=\infty$ unless  $t^2$ is a primitive $n$-th root of unity and $n>1$, in which case, $n_t=n$.

Note that if $n_t=n<\infty$, then $[n]_t=0$ and $[n-k]_t=-t^n[k]_t$ for $0\le k\le n$. In turn, this implies ${m\brack k}_t=0$ 
for all $m\ge n_t$, $1\le k\le m-1$ and 
$$
{n-1\brack k}_t=-t^n{n-1\brack k-1}_t
$$
hence 
$$
{n-1\brack k}_t=(-t^n)^k=1,
$$ 
which  most of the structure constants in  \eqref{eq:k-th divided power1} and \eqref{eq:k-th divided power2}. 
In particular, the following relations hold in $\tilde A_t$
$$
\sigma_1^{[k]}\sigma_1^{[n-k]}=\sigma_2^{[k]}\sigma_2^{[n-k]}=0, \quad \sigma_1^{[k]}\sigma_2^{[n-k]}=\sigma_{12}^{[n]}\ , $$
for all  $1\le k<n=n_t$,
where
$$
\sigma_{12}^{[n]}:=(-t^n)^{k-1}\sigma_1^{[n]}+(-t^n)^k\sigma_2^{[n]}\ .
$$

Now define the algebra $A_t$, $t\in \C^\times$ over $\Z[t+t^{-1}]\subset \C$ as follows:

If $n_t=\infty$, then $A_t:=\tilde A_t$;

If $n_t=n<\infty$ ( i.e., $t^2\ne 1$ is the $n$-th primitive root of unity), then $A_t$ is a subalgebra of $\tilde A_t$ generated by all 
$\sigma_1^{[k]},\sigma_2^{[k]}$, $k=0,1,\ldots,n-1$ and by $\sigma_{12}^{[n]}$.

It is easy to see that in both cases the algebra  $A_t$ is $\Z$-graded via $\deg \sigma_i^{[k]}=k$. Moreover,
in the second case,  $\deg \sigma_{12}^{[n]}=n$ is the top degree in $A_t$, as $[n]_t=0$.

For $t\in \C^\times$ let $W_t:=\langle s_1,s_2: s_1^2=s_2^2=1, (s_1s_2)^{n_t}=1\rangle$ be the dihedral group. 
Here it is understood that for $t=\pm 1$ we have the relation $s_1 s_2=1$ and for $t$ which is not a root of unity, we 
have the tautological relator  $(s_1s_2)^0=1$. Define the $W_t$-action on the {\em weight lattice}
$$
\La_t=\Z[t+t^{-1}]\cdot \sigma_1+\Z[t+t^{-1}]\cdot \sigma_2$$
by:
\begin{equation}
\label{eq:Wt action fundamental weights}
s_i(\sigma_j)=\sigma_j-\delta_{ij}(2\sigma_j-(t+t^{-1})\sigma_{3-i})
\end{equation}
for all $i,j\in \{1,2\}$.

Recall that if $W$ is a group acting on a vector space $V$, then the {\em coinvariant algebra}
$S(V)_W$ is the quotient $S(V)/\langle S(V)^W_+\rangle$, where $S(V)_+^W$ stands for all $W$-invariants
in the algebra of the constant-term-free polynomials  $S(V)_+=\sum_{k>0} S^k(V)$. (The computations of $S(V)_W$ 
below, in the case of $W_t$ with $t$ a root of unity, present a very special case of the computation of coinvariant 
algebras for arbitrary finite groups, see e.g. \cite{Hiller}.) 

The following proposition explains the origin of the algebra $A_t$:

\begin{proposition} For each $t\in \C^\times$ the algebra $\C\otimes A_t$ is naturally isomorphic to the
coinvariant algebra of $W_t$ acting on the vector space $V=\C \otimes \La_t$.
In particular, $W_t$ naturally acts on $A_t$
 via:
\begin{equation}
\label{eq:Wt action powers fundamental weights}
s_i(\sigma_j^{[k]})=\sigma_j^{[k]}-\delta_{ij}(2\sigma_j^{[k]}-(t^k+t^{-k})\sigma_{3-i}^{[k]})
\end{equation}
for all $i,j\in \{1,2\}$, $0\le k<n_t$ and (whenever $1<n_t=n<\infty$)
$$s_i(\sigma_{12}^{[n]})=-\sigma_{12}^{[n]}$$
for $i=1,2$.

\end{proposition}
\begin{proof} Denote $z_1=\sigma_1-t\sigma_2$, $z_2=t^{-1}\sigma_2-\sigma_1$ and let
$$e_2=-z_1z_2=(\sigma_1-t\sigma_2)(\sigma_1-t^{-1}\sigma_2)=
 \sigma_1^2+\sigma_2^2-(t+t^{-1})\sigma_1\sigma_2$$
(see \eqref{eq:quadratic coinvariant}). It is easy to see that under the action \eqref{eq:Wt action fundamental weights}, one has
\begin{equation}
\label{action}
s_1(z_1)=z_2,~s_1(z_2)=z_1,~s_2(z_1)=t^2z_2,~s_2(z_2)=t^{-2}z_2 \ .
\end{equation}
Hence, $e_2$ is invariant under the $W_t$-action.

Now assume that $[k]_t!\ne 0$ for all $k$, i.e.,  $n_t=\infty$. Then  the algebra $\C\otimes A_t$
is just the quotient of $\C[\sigma_1,\sigma_2]$ by the quadratic ideal generated by $e_2$.

On the other hand, it is easy to see, using \eqref{action},  that the $W_t$-invariant algebra $\C[\sigma_1,\sigma_2]^{W_t}$ is generated by $e_2$.
Therefore, the coinvariant algebra $\C[\sigma_1,\sigma_2]_{W_t}$ is also the quotient $\C[\sigma_1,\sigma_2]/\langle e_2\rangle$. This proves the proposition in the case when $n_t=\infty$.

Assume that now $n_t=n<\infty$ or, equivalently, $[k]_t!\ne 0$ for $k<n$ and $[k]_t!=0$ for $k\ge n$. Therefore, Proposition \ref{pr:k-th divided power1} guarantees that $\C\otimes A_t$ is a commutative algebra generated by $\sigma_1,\sigma_2$ subject to the relations
\begin{equation}\label{relators}
e_2=0,~\sigma_1^n=\sigma_2^n=0
\end{equation}
(In fact, $\sigma_{12}^{n}=\sigma_1\sigma_2^{n-1}$ by \eqref{eq:k-th divided power2} because $[n-1]_t=-t^n$.)

Again, it is easy to see, using \eqref{action}, that the $W_t$-invariant algebra $\C[\sigma_1,\sigma_2]^{W_t}$ is generated by $e_2$ and $e_n=z_1^n+z_2^n$.  Therefore, the coinvariant algebra $\C[\sigma_1,\sigma_2]_{W_t}$ is the quotient $\C[\sigma_1,\sigma_2]/\langle e_2,e_n\rangle$. To finish the proof it suffices to show that the ideals $\langle \sigma_1^n,\sigma_2^n\rangle$ and $\langle e_n\rangle$ are equal in $\C[\sigma_1,\sigma_2]/\langle e_2\rangle$.
Indeed, taking into account that $z_1z_2=0$ in $\C[\sigma_1,\sigma_2]/\langle e_2\rangle$ and that
$$\sigma_2=\frac{z_1+z_2}{t^{-1}-t},~\sigma_1=\frac{t^{-1}z_1+tz_2}{t^{-1}-t} \ ,$$
we obtain:
$$\sigma_1^n=\frac{(z_1+z_2)^n}{(t^{-1}-t)^n}=\frac{z_1^n+z_2^n}{(t^{-1}-t)^n},~\sigma_2^n=\frac{(t^{-1}z_1+tz_2)^n}{(t^{-1}-t)^n}=\frac{t^{-n}z_1^n+t^nz_2^n}{(t^{-1}-t)^n}=t^n\frac{z_1^n+z_2^n}{(t^{-1}-t)^n}$$
because $t^{2n}=1$. This proves the equality of ideals hence the equality of quotients $\C\otimes A_t=\C[\sigma_1,\sigma_2]_{W_t}$.

In particular, this verifies that $W_t$ naturally acts on $\C\otimes A_t$. To obtain \eqref{eq:Wt action powers fundamental weights}, note that $W_t$ preserves the component $\C\cdot \sigma_1^k+\C\cdot \sigma_2^k\subset \C\otimes A_t$ which is the $k$-th symmetric power of
$\C\cdot \sigma_1+\C\cdot \sigma_2$ (if $n<n_t$) and therefore, $W_t$ acts on the former space in the same way as in the latter space, i.e., by \eqref{eq:Wt action fundamental weights} where $t$ is replaced with $t^k$.

The proposition is proved.
\end{proof}

It is convenient to label the above basis of $A_t$ by the elements of the dihedral group $W_t$:
\begin{equation}
\label{eq:Schubert basis of A_t}
\sigma_w=
\begin{cases} \sigma_i^{[k]} &\text{if $\ell(w)=k<n_t$ and $\ell(ws_i)<\ell(w)$}\\
\sigma_{12}^{[n_t]} &\text{ if  $\ell(w)=n_t<\infty$}\\
\end{cases}
\end{equation}
for $w\in W_t$, where $\ell:W\to \Z_{\ge 0}$ is the length function.

The following result is an  equivalent reformulation of Proposition \ref{pr:k-th divided power1}.
\begin{proposition} For each $t\in \C^\times$ the elements $\sigma_w$, $w\in W_t$ form a $\Z[t+t^{-1}]$-linear basis of $A_t$ and the following relations are defining:

$\bullet$ If $\ell(u)+\ell(v)>n_t$, then $\sigma_u\sigma_v=0$.

$\bullet$ If $u=\underbrace{\cdots s_js_i}_k$, $v=\underbrace{\cdots s_js_i}_\ell$ and $k+\ell\le n_t$ and $\{i,j\}=\{1,2\}$, then
$$\sigma_u\sigma_v={k+\ell\brack k}_t\sigma_w$$
where $w={\underbrace{\cdots s_js_i}_{k+\ell}}$ (e.g., the right hand side is $0$ if $k+\ell=n_t$ and $k,\ell>0$).

$\bullet$ If $u=\underbrace{\cdots s_2s_1}_k$, $v=\underbrace{\cdots s_1s_2}_\ell$ and $k+\ell< n_t$, then
$$\sigma_u\sigma_v={k+\ell-1\brack k-1}_t\sigma_{w_1}+{k+\ell-1\brack \ell-1}_t\sigma_{w_2}$$
where
$$w_1={\underbrace{\cdots s_2s_1}_{k+\ell}},~w_2={\underbrace{\cdots s_1s_2}_{k+\ell}}$$.

$\bullet$ If $u=\underbrace{\cdots s_2s_1}_k$, $v=\underbrace{\cdots s_1s_2}_\ell$ and $k+\ell= n_t$, $k\le \ell$, then
$$\sigma_u\sigma_v={n_t-1\brack k-1}_t\sigma_{w_\circ}=\sigma_{w_\circ}$$
where $w_\circ={\underbrace{\cdots s_2s_1}_{n_t}}$ is the longest element of the (finite) group $W_t$.

\end{proposition}

%
%
%

Note that when $\theta=t+t^{-1}\in \R$,  all structure constants of $A_t$ are real numbers. We can refine this as follows.

\begin{corollary} 
The structure constants of $A_t$  are non-negative if and only if either $t=e^{\frac{\pi\sqrt{-1}}{n}}$ or $t>0$.
\end{corollary}
\begin{proof} Indeed, the structure constants are $t$-binomials, which are non-negative for $t=e^{\frac{\pi\sqrt{-1}}{n}}$ or $t>0$, 
since $[m]_t\ge 0$ for $1\le m \le n$. On the other hand, if, say, $n=n_t<\infty$ but $t$ is not of the form $e^{\frac{\pi\sqrt{-1}}{n}}$, then there exists 
$1\le m \le n$ so that $[m]_t<0$. \qed \end{proof}

\begin{rem}
The above corollary is just one of many hints pointing to existence of (possibly non-commutative, in view of non-integrality of the structure constants) 
complex-algebraic varieties serving as flag-manifolds for non-crystallographic finite dihedral groups. 
\end{rem}

Let $G$ be a complex Kac-Moody group  with  the Cartan matrix
$$
\begin{pmatrix}
2 & -a_{12}\\
-a_{21} & 2
\end{pmatrix},
$$
where $a_{12}$ and $a_{21}$ are arbitrary positive integers (if $a_{12}a_{21}\le 3$, then $G$ is a finite-dimensional simple Lie group of rank $2$). Let $t\in \C^\times$ be
such that $t+t^{-1}=\sqrt{a_{12}a_{21}}$. In particular, the Weyl group of $G$ is naturally isomorphic to $W_t$.
It is well-known (see e.g., \cite{koku})
that the cohomology algebra $H^*(G/B)$ has a basis of  Schubert classes $[X_w]$, $w\in W_t$.

The following is the main result of the section.

\begin{theorem}
\label{th:cohomology} Let $G$ and $B$ be as above and $c_1,c_2\in \C^\times$ be any numbers such that $\frac{c_1}{c_2}=\sqrt{\frac{a_{12}}{a_{21}}}$ and 
$\Z[c_1,c_2]\supset \Z[t+t^{-1}]$. Then the association
\begin{equation}
\label{eq:the isomorphism}
[X_w]\mapsto c_i^{\lceil \frac{k}{2}\rceil}c_{3-i}^{\lfloor \frac{k}{2}\rfloor}\cdot \sigma_w
\end{equation}
for all $w\in W_t$, where $i\in \{1,2\}$ is such that $\ell(ws_i)<\ell(w)=k$, defines a  $W_t$-equivariant isomorphism
\begin{equation}
\label{eq:the canonical isomorphism}
H^*(G/B,\Z[c_1,c_2])\widetilde \to\Z[c_1,c_2]\otimes A_t \ .
\end{equation}
\end{theorem}

\begin{proof} It suffices to prove that \eqref{eq:the isomorphism} defines a $W_t$ equivariant isomorphism
\begin{equation}
\label{eq:isomorphism cohomology}
H^*(G/B,\C)\widetilde \to \C\otimes A_t
\end{equation}

Recall that the action of the Weyl group $W$ of $G$ on the root space $Q_{\C}=\C\cdot \alpha_1+\C\cdot \alpha_2$ is given by:
$$s_i(\alpha_j)=\begin{cases}
-\alpha_i & \text{if $i=j$}\\
\alpha_i+a_{ij}\cdot  \alpha_j & \text{if $i\ne j$}\\
\end{cases}$$
for $i,j\in \{1,2\}$.

It follows from  \cite[Proposition 3.10]{koku} that the algebra $H^*(G/B,\Z)$  satisfies the following Chevalley formula:
\begin{equation}
\label{eq:chevalley}
 [X_w][X_{s_i}] =\sum_{w_1,w_2,j} \omega_i^\vee(w_2^{-1}(\alpha_j))\cdot [X_{w_1s_jw_2}] \ ,
\end{equation}
where the summation is over all $w_1,w_2\in W$, and $j\in \{1,2\}$ such that $w=w_1w_2$, $\ell(w)=\ell(w_1)+\ell(w_2)$, and $\ell(w_1s_jw_2)=\ell(w)+1$. Here $\omega_i^\vee$, $i\in \{1,2\}$, denotes the dual basis in $Q^*$ of the basis
$\alpha_1,\alpha_2$. 

In particular, if $\ell(ws_i)<\ell(w)$, then the only non-zero summand in the right hand side of \eqref{eq:chevalley} corresponds to $w_1=1$ and $w_2=w$, and $j$ such that $\ell(s_jw)=\ell(w)+1$. Furthermore,  if $\ell(ws_i)>\ell(w)$, then the right hand side has two summands, first of which comes with $w_2=1$, $w_1=w$, and the second one -- with $w_1=1$, $w_2=w$.
Therefore,
\begin{equation}
\label{eq:Chevalley dihedral1}
[X_{\underbrace{\cdots s_js_i}_k}][X_{s_i}]=\omega_i^\vee(\underbrace{s_is_j\cdots s_{i'}}_k(\alpha_{3-i'}))[X_{\underbrace{\cdots s_js_i}_{k+1}}]\ ,
\end{equation}
\begin{equation}
\label{eq:Chevalley dihedral2}
[X_{\underbrace{\cdots s_is_j}_k}][X_{s_i}]=[X_{\underbrace{\cdots s_js_i}_{k+1}}]+\omega_i^\vee(\underbrace{s_is_j\cdots s_{i'}}_k(\alpha_{3-i'}))[X_{\underbrace{\cdots s_is_j}_{k+1}}]
\end{equation}
for all $k<n_t$ and $i,j$ such that $\{i,j\}=\{1,2\}$, where $i'$ stands for the appropriate index $i$ or $j$ (depending on $k \mod 2$).
In particular, if $k=1$, we obtain:
$$[X_{s_1}]^2=\omega_1^\vee(s_1(\alpha_2))[X_{s_2s_1}],~[X_{s_2}]^2=\omega_2^\vee(s_2(\alpha_1))[X_{s_2s_1}],[X_{s_1}][X_{s_2}]=[X_{s_1s_2}]+[X_{s_2s_1}] $$
which implies the following quadratic relation in $H^*(G/B,\Z)$:
\begin{equation}
\label{eq:quadratic relation}
a_{21}[X_{s_1}]^2+a_{12}[X_{s_2}]^2=a_{12}a_{21}[X_{s_1}][X_{s_2}]
\end{equation}

To utilize the identities \eqref{eq:Chevalley dihedral1} and \eqref{eq:Chevalley dihedral2}, we need the following obvious result.

\begin{lemma}
\label{le:action dihedral}
Let $w=\underbrace{\cdots s_is_j}_k\in W_t$, where $\{i,j\}=\{1,2\}$. Then
$$w(\alpha_i)=[k+1-\varepsilon_k]_t\alpha_i +\sqrt{\frac{a_{ji}}{a_{ij}}}\cdot [k+\varepsilon_k]_t\alpha_j \ ,$$
where $t+t^{-1}=\sqrt{a_{12}a_{21}}$ and $\varepsilon_k=\begin{cases}
1 & \text{if $k$  is odd}\\
0 & \text{if $k$  is even}\\
\end{cases}$.
\end{lemma}

Therefore, we obtain:
\begin{equation}
\label{eq:Chevalley dihedral special1}
[X_{\underbrace{\cdots s_js_i}_k}][X_{s_i}]=\left(\sqrt{\frac{a_{ij}}{a_{ji}}}\right)^{\varepsilon_k}[k+1]_t\cdot [X_{\underbrace{\cdots s_js_i}_{k+1}}]
\end{equation}
\begin{equation}
\label{eq:Chevalley dihedral special2}
[X_{\underbrace{\cdots s_is_j}_k}][X_{s_i}]=[X_{\underbrace{\cdots s_js_i}_{k+1}}]+\left(\sqrt{\frac{a_{ij}}{a_{ji}}}\right)^{\varepsilon_k}[k+1]_t\cdot [X_{\underbrace{\cdots s_is_j}_{k+1}}]
\end{equation}
Furthermore, \eqref{eq:Chevalley dihedral special1} implies that
\begin{equation}
\label{eq:Xi to k}
[X_{s_i}]^k=\left(\sqrt{\frac{a_{ij}}{a_{ji}}}\right)^{\lfloor \frac{k}{2}\rfloor}[k]_t!\cdot [X_{\underbrace{\cdots s_js_i}_k}] \ .
\end{equation}
In turn, this implies that $H^*(G/B,\C)$ is generated by $[X_{s_1}]$, $[X_{s_2}]$  satisfying \eqref{eq:quadratic relation} and the relations
\begin{equation}\label{finiteorder}
[X_{s_1}]^{n_t}=[X_{s_2}]^{n_t}=0\end{equation}
if  $n_t<\infty$.
Pick $r_1,r_2\in \C^\times$ such that $\frac{r_1}{r_2}=\sqrt{\frac{a_{21}}{a_{12}}}$ and define
$$\varphi:\sigma_1\mapsto r_1 [X_{s_1}],\quad \varphi:\sigma_2\mapsto r_2 [X_{s_2}].$$
In view of the relation \eqref{eq:quadratic relation}, we obtain
$$
\varphi(\si_1)^2 + \varphi(\si_2)^2 = \sqrt{a_{12} a_{21}} \varphi(\si_1\si_2).
$$
Since $t+t^{-1}= \sqrt{a_{12} a_{21}}$, we conclude that $\varphi$ preserves the
defining quadratic equation \eqref{eq:quadratic coinvariant} of $A_t$. The equation
\eqref{finiteorder} implies that $\varphi$ preserves the last two relators in \eqref{relators} provided that $n=n_t<\infty$. Thus, $\varphi$ extends a surjective homomorphism of algebras $\varphi:\C\otimes A_t\to H^*(G/B,\C)$.

Clearly, this homomorphism is an
isomorphism because it preserves the natural $\Z$-grading and the respective graded components of both algebras are of the same dimension. Furthermore, let us show that
for each $w\in W_t$ one has:
\begin{equation}
\label{eq:image Schubert classes}
\varphi(\sigma_w)=r_i^{\lceil \frac{\ell(w)}{2}\rceil}r_j^{\lfloor \frac{\ell(w)}{2}\rfloor}[X_w]^{\ell(w)}
\end{equation}
where $i\in \{1,2\}$ is such that $\ell(ws_i)<\ell(w)$ and $\{i,j\}=\{1,2\}$. Indeed, if $\ell(w)<n_t$, then
$$\varphi(\sigma_w)=r_i^{\ell(w)}\frac{1}{[k]_t!}[X_{s_i}]^{\ell(w)}=r_i^{\ell(w)}\left(\frac{r_j}{r_i}\right)^{\lfloor \frac{\ell(w)}{2}\rfloor}[X_w]^{\ell(w)}=r_i^{\lceil \frac{\ell(w)}{2}\rceil}r_j^{\lfloor \frac{\ell(w)}{2}\rfloor}[X_w]^{\ell(w)} \ .$$
If  $\ell(w)=n_t<\infty$ (i.e., $w$ is the longest element of $W$), then $[n_t]_t=0$ and, using \eqref{eq:k-th divided power2} with $k=1$ and \eqref{eq:Chevalley dihedral special2} respectively, we obtain:
$$\sigma_{ws_1}\sigma_1=\sigma_w,~[X_{ws_1}][X_{s_1}]=[X_w] \ .$$
Thus applying $\varphi$ to the first of these relations, we obtain (taking into account that $r_1=r_2$ when $n_t$ is odd and using the already proved case of \eqref{eq:image Schubert classes} with $w'=ws_1$, $i=2$):
$$\varphi(\sigma_w)=\varphi(\sigma_{ws_1}\sigma_1)=\varphi(\sigma_{ws_1})\varphi(\sigma_1)=r_2^{\lceil \frac{n_t-1}{2}\rceil}r_1^{\lfloor \frac{n_t-1}{2}\rfloor}[X_{ws_1}]r_1[X_{s_1}]=(r_1r_2)^{ \frac{n_t}{2}}[X_{w}] \ .$$
Finally, taking $r_i=\frac{1}{c_i}$, $i=1,2$, we see that the isomorphism $\varphi^{-1}$ is given by \eqref{eq:the isomorphism} and its restriction to $H^*(G/B,\Z[c_1,c_2])$ becomes \eqref{eq:the canonical isomorphism}. The $W$-equivariancy of both $\varphi$ and $\varphi^{-1}$ follows.

\end{proof}

\begin{rem}
A computation of the rings $H^*(G/B, \Z)$ for rank 2 complex Kac-Moody groups $G$ appeared in 
\cite[Section 10]{Kit}. We are grateful to Shrawan Kumar for this reference. 
\end{rem}

\begin{remark} We can take
$$
c_i=\sqrt{\frac{a_{i,3-i}}{gcd(a_{12},a_{21})}}, \quad i=1,2$$
in Theorem \ref{th:cohomology}. Then $\Z[c_1,c_2]\supset \Z[t+t^{-1}]$ because $t+t^{-1}=c_1c_2 \cdot gcd(a_{12},a_{21})$.
In particular, if the Cartan matrix is symmetric, i.e., $a_{12}=a_{21}$, then the isomorphism \eqref{eq:the canonical isomorphism} is over $\Z$  because $c_1=c_2=1$ and $\Z[c_1,c_2]=\Z[t+t^{-1}]=\Z$.
\end{remark}

In view of Theorem \ref{th:cohomology}, we will refer to the algebra $A_t$ as the {\it universal dihedral cohomology} and to 
the basis $\{\sigma_w\}$ -- as the {\it universal Schubert classes}: Under various specializations of $t$ it computes either cohomology rings of complex 
flag manifolds associated with complex Kac-Moody groups or cohomology rings of ``yet to be defined'' flag-manifolds for non-crystallographic finite dihedral groups or nondiscrete 
infinite dihedral groups. 

We call a complex number $t$ {\it admissible} if  either

(1) (finite case) $t=e^{\pm \frac{\pi\sqrt{-1}}{n}}$ for some $n\in \Z_{>0}$, or

(2) (hyperbolic case) $t$ is a positive real number.

\medskip
Then for every admissible $t$, $[k]_t> 0$ for all $0\le k<n_t$. For an admissible $t$ and let $W_t^{(i)}=\{w\in W_t|\ell(ws_{3-i})=\ell(w)+1\}$, $i=1, 2$. 

\begin{notation}
Denote by  $B_t^{(i)}, i=1, 2,$ be the subalgebras of $A_t$ generated by $X^{(i)}=\{\si_w: w\in W_t^{(i)}\}$. 
\end{notation}

The subalgebras $B_t^{(i)}$ play the role of the cohomology rings of the ``Grassmannians'' $Y_i, i=1, 2$.  
It follows from Proposition \ref{pr:k-th divided power1} that $X^{(i)}$ is a basis of 
$B_t^{(i)}$, e.g., $\dim B_t^{(i)}=|W_t^{(i)}|=n_t-1$, and that, moreover, the ring $B_t^{(i)}$ is naturally isomorphic to the cohomology ring
$$
H^*(\C \P^n), \quad n=n_t. 
$$

\section{Belkale-Kumar type filtration of $A_t$}\label{sec:BKfiltrations}

In this section, we construct a filtration on $A_t$ (and it subalgebras $B_t^{(i)}, i=1, 2$) in the sense of Proposition
\ref{pr:crystal limit Belkale-Kumar additive} using a Belkale-Kumar type function $\varphi:W_t\to \R$.
In the case, when $t$ is $n$-th primitive root of unity, the associated graded algebra $A_{t,0}=gr A_t$ will play a role of Belkale-Kumar 
cohomology of spherical buildings with finite Weyl group $I_2(n)$,  
which is ``Poincar\'e dual'' to the homology pre-ring $H_*(X)$ defined by the Schubert pre-calculus.

\begin{definition}
\label{def:concave}
Let $\k$ be a field and $A$ be an associative $\k$-algebra with a basis $\{b_x|x\in X\}$ so that
\begin{equation}
\label{eq:associative algebra}
b_x b_y= \sum_{z\in X} c^z_{x,y} b_z.
\end{equation}
for all $x,y\in X$, where $c_{x,y}^z\in \k$ are structure constants. Furthermore, given an ordered abelian semi-group $\Gamma$
(e.g., $\Gamma=\R$), we say that a function $\varphi:X\to \Gamma$ is {\it concave} if
\begin{equation*}
\varphi(x)+\varphi(y)\ge \varphi(z)
\end{equation*}
for all $x,y,z\in X$ such that $c_{x,y}^z\ne 0$.
\end{definition}

\begin{proposition}
\label{pr:crystal limit Belkale-Kumar additive} In the notation \eqref{eq:associative algebra}, for each concave function $X\to \Gamma$ we have:

(a)  $A$ is filtered by $\Gamma$ via $A_{\le \gamma}:=\sum\limits_{x\in X:\varphi(x)\le \gamma} \k\cdot b_x$.

(b) The multiplication in the associated graded algebra $A_0=gr A$ is given by:
\begin{equation}
\label{eq:crystal limit Beklale-Kumar additive}
b_x\circ b_y= \sum_{z\in X:\varphi(z)=\varphi(x)+\varphi(y)} c^z_{xy} b_z.
\end{equation}
for all $x,y\in X$, where $c_{xy}^z\in \k$ are the structure constants of $A$.

\end{proposition}

\begin{proof} Part (a). Assume that $\varphi(x)\le \gamma_1$, $\varphi(y)\le \gamma_2$, i.e.,
$b_x\in A_{\le \gamma_1}$, $b_y\in A_{\le \gamma_2}$. Then each $z$
such that $c_{xy}^z\ne 0$ satisfies $\varphi(z)\le \varphi(x)+\varphi(y)\le \gamma_1+\gamma_2$,
i.e., $b_z\in A_{\gamma_1+\gamma_2}$. Therefore, $b_xb_y\in A_{\gamma_1+\gamma_2}$. This proves (a).
Part (b) immediately follows. \qed 
\end{proof}

\begin{rem}\label{bk:rem}
The algebra $A_0$ is the {\em Belkale-Kumar degeneration} of $A$. It was introduced by Belkale and Kumar in \cite{bk} in the special case 
of cohomology rings of flag-manifolds $G/B$, with $G$ a complex semisimple Lie group and $B$ its Borel subgroup. In order to relate our
definition to that of \cite{bk}, note that, given a concave function $\varphi$, Belkale and Kumar define the deformation $A_\tau$ of 
$A=H^*(G/B,\C)$ by
$$
b_x\odot_\tau b_y:= \sum_{z\in X} \tau^{\varphi(x)+\varphi(y)-\varphi(z)}c^z_{xy} b_z. 
$$
Setting $\tau=1$ one recovers the original algebra $A$, while sending $\tau$ to zero one obtains the degeneration $A_0=gr(A)$ of $A$.   
\end{rem}

Our goal is to generalize the function $\varphi$ defined in \cite{bk} to the case of algebras $A_t$ (for {\em admissible} values of $t$), 
so that our function $\varphi$ will specialize to the Belkale-Kumar function in the case $n=3, 4, 6$. 
Note that concavity of $\varphi$ was proven in \cite{bk} as a consequence of complex-algebraic nature of the variety $G/B$. 
In our case, such variety does not exist and we prove concavity 
by a direct calculation.    

\medskip
For $t\in \C^\times$ define the action of the dihedral group $W_t$ on the $2$-dimensional {\em root lattice}
$$
Q=Q_t=\Z[t+t^{-1}]\cdot \alpha_1+ \Z[t+t^{-1}]\cdot \alpha_2$$
by:
$$s_i(\alpha_j)=\begin{cases}
-\al_i & \text{if $i=j$}\\
\al_i+[2]_t\cdot  \al_j & \text{if $i\ne j$}\\
\end{cases}$$
for $i,j\in \{1,2\}$.

For each $i\in \{1,2\}$ define the map $[\cdot]_i:W_t\to Q$ recursively by $[1]_i=0$ and
$$[s_jw]_i=\delta_{ij}\alpha_i+s_j([w]_i) \ .$$

The above action extends to the {\em weight lattice}
$$
\Lambda=\Lambda_t=\Z[t+t^{-1}]\cdot \omega_1+ \Z[t+t^{-1}]\cdot \omega_2$$
by:
$$
s_i(\omega_j)=\omega_j-\delta_{ij}\iota(\alpha_i)$$
for all $i,j\in \{1,2\}$, which is consistent with \eqref{eq:Wt action fundamental weights}. 
Here $\iota:Q\to \Lambda$ is a $\Z[t+t^{-1}]$-linear map given by:
$$
\iota(\alpha_1)=2\omega_1-(t+t^{-1})\omega_2,~\iota(\alpha_2)=2\omega_2-(t+t^{-1})\omega_1 \ .
$$
Note that the above map $[\cdot]_i$ satisfies:
$$\iota([w]_i)=\omega_i-w(\omega_i) \ .$$
Define the functions $\Phi_i:W_t\to \Z[t+t^{-1}]$, $i=1,2$ by
\begin{equation}
\label{eq:Belkale-Kumar convex function}
\Phi_i(w)=|[w]_i| 
\end{equation}
where $|g_1 \alpha_1+g_2 \alpha_2|=g_1+g_2$.
\begin{proposition} 
\label{pr:Phi_i}
For any $w\in W_t$, $i=1,2$ we have:
\begin{equation}
\label{eq:Phi_i}
\Phi_i(w)=\begin{cases}
{\ell(w)+1\brack 2}_q &\text{if $\ell(ws_i)<\ell(w)$}\\
{\ell(w)\brack 2}_q &\text{if $\ell(ws_i)>\ell(w)$}\\
\end{cases},
\end{equation}
where $q=t^{1/2}$ and $\ell: W\to \Z$ is the word-length function on $W$ with respect to the generating set $s_1, s_2$. 
 
In particular, the function $\Phi:=\Phi_1+\Phi_2$  is given by the formula:
\begin{equation}
\label{eq:Phi total}\Phi(w)=([\ell(w)]_q)^2.
\end{equation}

\end{proposition}

\begin{proof} We need the following obvious result:

\begin{lemma}
\label{le:rho}
For each $k\in \Z$ denote
$\alpha_k:=\begin{cases} \alpha_1 & \text{if $k$ is odd}\\
\alpha_2 & \text{if $k$ is even}\\
\end{cases} .$

Let $w=\underbrace{\cdots s_js_i}_k\in W_t$, where $\{i,j\}=\{1,2\}$. Then
$$w(\alpha_j)=[k]_t\alpha_{i+k}+[k+1]_t\alpha_{j+k},~[w]_i=\alpha_i+[2]_t\alpha_{i+1}+\cdots +[k]_t\alpha_{i+k-1}.$$
\end{lemma}
\begin{proof} The assertion directly follows from Lemma \ref{le:action dihedral} with $a_{12}=a_{21}=t+t^{-1}$. \qed 
\end{proof}

Furthermore, using the second identity of Lemma \ref{le:rho} we obtain for any $w\in W_t$ with $\ell(ws_i)<\ell(w)$:
$$|[w]_i|=[1]_t+[2]_t+\cdots+[\ell(w)]_t\ .$$
Using the fact that $[m]_t=\frac{[2m]_q}{[2]_q}$ for $q=t^{1/2}$ and any $m$, we obtain
$$|[w]_i|=[1]_t+[2]_t+\cdots+[\ell(w)]_t=\frac{1}{[2]_q}([2]_q+[4]_q+\cdots+[2\ell(w)]_q)=\frac{1}{[2]_q}[\ell(w)]_q[\ell(w)+1]_q\ , $$
which  proves \eqref{eq:Phi_i} since $\Phi_i(w)=\Phi_i(ws_{3-i})=|[w]_i|$.

We now prove \eqref{eq:Phi total}. Indeed, for any $w\in W_t$ let $i$ be such that $\ell(ws_i)<\ell(w)$. Applying part \eqref{eq:Phi_i}, we obtain:
$$\Phi(w)=|[w]_i|+|[w]_{3-i}|={\ell(w)+1\brack 2}_q+{\ell(w)\brack 2}_q=([\ell(w)]_q)^2\ .$$
The proposition is proved. \qed 
\end{proof}

The following theorem is the main result of the section.

\begin{theorem}  \label{thm:filtrations}
The functions
$$
\varphi_i: X^{(i)}\to \R, \quad \varphi_i(\si_w)= -\Phi_i(w)
$$
$i=1,2$ and 
$$
\varphi: X\to \R, \quad \varphi(\si_w)= -\Phi(w)
$$
(see Proposition \ref{pr:Phi_i}) are both concave in the sense of Definition \ref{def:concave}; in particular, they define filtrations on $B_t^{(i)}$ and $A_t$ respectively in the sense of Proposition \ref{pr:crystal limit Belkale-Kumar additive}.
Moreover, the equalities
$$
\varphi(\si_u)+\varphi(\si_v)=\varphi(\si_w), \quad
\varphi_i(\si_u)+\varphi_i(\si_v)=\varphi_i(b_w)
$$
are achieved if and only if  either:

1. For the function $\varphi$, $u=1$ or $v=1$, or $\ell(u)+\ell(v)=\ell(w)=n$, provided that $n<\infty$.

2. For the function $\varphi_i$, $u=1$ or $v=1$, or $\ell(u)+\ell(v)=\ell(w)=n-1$, provided that $n<\infty$.

\end{theorem}
\proof Recall that a function $f: I \to \R$ defined on an interval $I\subset \R$ is called {\em superadditive} (resp. subadditive) if
\begin{equation}\label{eq:super}
f(x+y)\ge f(x)+ f(y), \quad \hbox{resp.}\quad f(x+y)\le f(x)+ f(y)
\end{equation}
for all $x, y, x+y\in I$. If $f$ is convex, continuous, and $f(0)=0$ then $f$ is superadditive on
$I=\R_+$, see \cite[Theorem 7.2.5]{HP}. Moreover, it follows from the proof of
\cite[Theorem 7.2.5]{HP} that if $f$ is strictly convex then \eqref{eq:super} is a strict inequality unless $xy=0$.

Let $t\in \C$ be an admissible number, $n:=n_t$; let $q:=t^{1/2}$ so that  $q\in \R_+$ if $t>0$ and
$q=e^{\sqrt{-1}Q}$, $Q=\frac{\pi}{2n}$ if $t$ is a root of unity. Define the functions
$$
F(x)= { x+1 \brack 2}_q, \quad 0\le x\le n-1,
$$
$$
G(x):= ([x]_q)^2, \quad 0\le x\le n,
$$
where $x$ are nonnegative real numbers. 

\begin{prop}\label{super}
The functions $F$ and $G$ are superadditive. Moreover, the inequality \eqref{eq:super} 
is  equality iff $xy(n-1-x-y)=0$ (for the function $F$) and $xy(n-x-y)=0$ (for the function $G$).
\end{prop}
\proof We have
$$
F(x)= \frac{[x]_q [x+1]_q}{[2]_q}= \frac{f(x)}{(q-q^{-1})^2(q+q^{-1})}, \quad f(x)=(q^x- q^{-x})(q^{x+1}- q^{-x-1})
$$
$$
G(x)= \frac{g(x)}{(q- q^{-1})^2}, \quad g(x)= (q^x- q^{-x})^2.
$$
In particular, $F(0)=G(0)=0$.

We first consider the hyperbolic case (i.e., $q>0$). Then the denominators of both $F$ and $G$ are positive and numerators are equal to
$$
f(x)=q^{2x+1}+q^{-2x-1}-q- q^{-1}
\quad g(x)= q^{2x}+ q^{-2x}-2
$$
It is elementary that both functions are strictly convex on $[0, \infty)$ because $f''(x)>0$ and $g''(x)>0$. Hence,  $F$ and $G$ are superadditive
with equality in \eqref{eq:super} iff $xy=0$.

\medskip
We therefore assume now that $q$ is a root of unity. One can check that in this case $F$ and $G$ are neither convex nor concave on
their domains, so we have to use a direct calculation in order to show superadditivity. The denominators of the functions $F$ and $G$
are both negative since they equal to
$-8 \sin^2(Q) \cos(Q)$ and $-4 \sin^2(Q)$ respectively.

Consider the functions $f(x)$ and $g(x)$. It is easy to see that 
$$
f(z)-f(x)-f(y)= (q^x- q^{-x})(q^y -q^{-y})(q^{x+y+1}+ q^{-x-y-1}),$$
$$g(z) -g(x)-g(y)= (q^x -q^{-x})(q^y - q^{-y})(q^z+ q^{-z}).$$
Therefore, if
$x, y,z \in [0, n-1]$ with $z=x+y$ then:
$$f(z) -f(x)-f(y)\le 0$$ with equality iff $xy(n-1-z)=0$ and 
$$g(z) -g(x)-g(y)\le 0$$ with equality iff $xy(n-z)=0$
because 
$$
(q^x- q^{-x})(q^y -q^{-y})=-4\sin(Qx)\sin(Qy)\le 0,~q^{u}+ q^{-u}=2\cos(Qu)\ge 0 
$$
for any $x,y,u\in [0,n]$.

\medskip
Thus both functions $f$ and $g$ are subadditive. Since the denominators in $F$ and $G$ are constant and negative, 
these functions are superadditive with equality in \eqref{super} iff $xy=0$ or $x+y=n-1$ (for $F$) and $x+y=n$ (for $G$). \qed

We can now finish the proof of Theorem \ref{thm:filtrations}. We have
$$
\Phi_i(w):=|[w]_i|=F(\ell(w)), \quad w\in W^{(i)}, 0\le \ell(w)\le n-1
$$
and
$$
\Phi(w)=|[w]_1|+ |[w]_2|= G(\ell(w)), w\in W, 0\le \ell(w)\le n.
$$
Observe that, since $A_t$ is graded by the length function of $W_t$,  
$$
c_{uv}^w\ne 0 \Rightarrow \ell(w)=\ell(u)+\ell(v),
$$
where $c_{uv}^w$ are the structure constants: 
$$
\si_u \cdot \si_v= \sum_{w} c_{uv}^w \si_w.  
$$
Therefore, superadditivity of the functions $F$ and $G$ is equivalent to concavity of the functions $\varphi=-\Phi$ and 
$\varphi_i=-\Phi_i$. The equality cases in Theorem \ref{thm:filtrations} immediately follow as well. \qed

\begin{cor}\label{bk-filtrations}
The rings $A_t$ and $B_t^{(i)}$, $i=1, 2$, admit Belkale--Kumar degenerations  $gr(A_t)$ and $gr(B_t^{(i)})$ 
given by the functions $\varphi$ and $\varphi_i$ respectively. 
\end{cor}

\begin{rem}
We do not know a natural topological interpretation for the rings  $gr(A_t)$ and $gr(B_t^{(i)})$. 
\end{rem}



\section{Interpolating between  homology pre-ring and the ring $gr(A_t)$}

Let $\k$ be a field. In this section we construct an interpolations between the  homology pre-rings 
$H_*(X, \widehat\k), H_*(X_l, \widehat\k)$ and the $\k$-algebras $gr(A_t), gr(B^{(l)}_t)$,  
$t=e^{\frac{\pi}{n}\sqrt{-1}}$, which are Belkale-Kumar degenerations of $A_t, B^{(l)}_t$ 
introduced in Section \ref{sec:BKfiltrations}. Thereby, we link the geometrically defined homology pre-rings and the 
algebraically defined cohomology rings of $X, X_l$, $l=1, 2$. 

Below we again abuse the terminology and use the notation $\infty$ for the infinity in the one-point compactification 
of $\R$ and for the element of $\widehat{\k}$. 
Accordingly, we equip $\k$ with the discrete topology and set 
$$
\lim_{\tau\to\infty} f(\tau) a=\infty, 
$$
whenever $a\in \k^\times$ and $\lim_{\tau\to\infty} f(\tau)=\infty$.

\medskip 
1. Interpolation for $A_t$. Using the Belkale-Kumar function $\varphi=-\Phi$ as in the previous section, we define the (trivial) 
family of algebras $A_{t,\tau}$ as in Remark \ref{bk:rem}, with  multiplication given (for $\tau>0$) by 
$$
\si_u\odot_\tau \si_v:= \sum_{w: \ell(w)=\ell(u)+\ell(v)} \tau^{\varphi(u)+\varphi(v)-\varphi(w)}c^w_{uv} \si_w \ ,  
$$
where $c_{uv}^w$ are the structure constants in $A_t$. Then, as $\tau\to 0$, the algebra $A_{t,\tau}$ degenerates to $gr(A_t)$. 
Now, let $\tau\to \infty$. Recall that $\varphi(u)+\varphi(v)-\varphi(w)>0$ unless it equals to zero (Proposition \ref{super}); 
the latter corresponds to the {\em degenerate cases}, i.e., 
products of Poincar\'e dual classes $\si_u, \si_v$ or classes where $\si_u=1$ or $\si_v=1$. 
Therefore, the limit pre-ring $A_{t,\infty}$ has structure constants $\hat{c}_{uv}^w$ equal to $0, 1, \infty$.

Here $\hat{c}_{uv}^w=0$ occurs unless $\ell(w)=\ell(u)+\ell(v)$, and $u, v, w\in W^{(i)}, i=1, 2$; 
in the latter case $\hat{c}_{uv}^w=\infty$ except for the degenerate 
cases where the structure constants are equal to $1$. Hence, in view of Proposition \ref{H_*(X)}, 
we obtain a degree-preserving isomorphism of  pre-rings  $A_{t,\infty}\cong H_*(X, \widehat{\k})$ given by
$$
\si_w \mapsto C_{w_\circ w}, \quad w\in W.  
$$ 

\medskip 
2. Interpolation for $B_t^{(l)}$, $l=1, 2$. The argument here is identical to the case of $A_t$, 
except the isomorphism is given by
\begin{equation}\label{BKPR}
\si_w\mapsto C_{n-1-r}\in H_*(X_l, \widehat{\k}), \quad r=\ell_l(w).    
\end{equation}

\medskip 
We conclude that the relation between $H_{BK}^*(X, \k):=gr(A_{t})$ and $H_*(X, \widehat{\k})$, 
is that of ``mirror partners'': They are different degenerations of a common ring $A_t$.

\section{Strong triangle inequalities}\label{sec:sti}

In this section we introduce a redundant system of inequalities
equivalent to $WTI$: This equivalence will be using in the following section.

Let $W=I_2(n)$ with the affine Weyl chamber $\De\subset \R^2$, $\k$ a field and $\widehat\k$ the corresponding pre-ring.  
Define the subset 
$$
\Si_{A,m}\subset W^m
$$
consisting of $m$-tuples $(u_1,...,u_m)$ of elements of $W$ so that
\begin{equation}\label{uprod}
\prod_{i} C_{u_i}=a\cdot C_{\1}, a\in \widehat{\k}^\times 
\end{equation}
in the pre-ring $H_*(X, \widehat{\k})$, where $X$ is a thick spherical building with the Weyl group $W$ satisfying
Axiom A. We then define  cones $K(\Si_{A,m})\subset \De^m$ by imposing the inequalities
$$
 \sum_{i} u_i^{-1}(\la_i) \le_{\De^*} 0 
$$
for the  $m$-tuples $(u_1,...,u_m)\in \Si_{A,m}$. 
We will refer to the defining inequalities of $K(\Si_{A,m})$ as {\em Strong Triangle Inequalities}, $STI$. 

Recall that $\K_m=\K_m(X)\subset \De^m$ is the stability cone of $X$, cut out by 
the inequalities $WTI$, see \S \ref{sec:sineq}. Then, clearly,
$$
K(\Si_{A,m})\subset  \K_m
$$
since the system $STI$ contains the $WTI$. The following is the main result of this section:

\begin{thm}\label{equalcones}
$$
K(\Si_{A,m})= \K_m.
$$
\end{thm}
\proof Observe that $u_i\ne \1$ for $i=1,...,m$, for otherwise the product in the left-hand side of \eqref{uprod} is zero.  
We first establish some inequalities concerning the relative lengths of elements of $W$: 

\begin{prop}\label{ine}
Let $w_i\in W\setminus \{\1\}, i=1,...,m$ are such that
$$
\prod_{i=1}^m C_{w_i}= C_\1 
$$ 
in the pre-ring $H_*(X, \widehat{\k})$.
Then for $k=1, 2$, we have:
$$
\sum_{i=1}^m \ell_k(w_i)\ge (m-1)(n-1). 
$$
In other words, for $r_i:=\ell_k(w_i)$, 
$$
\prod_{i=1}^m C_{r_i} \ne 0
$$
in the pre-ring $H_*(X_k, \widehat\k)$. 
\end{prop}
\proof Let $u_i, u\in W^{(j)}$ be such that 
$$
\prod_{i=1}^s C_{u_i}= a C_u, a\ne 0 
$$
in the pre-ring $H_*(X, \widehat\k)$. Then
$$
\sum_{i=1}^s (n-\ell(u_i))= n-\ell(u).  
$$
Since $\ell(u_i)=\ell_k(u_i)+ \del_{jk}, \ell(u)=\ell_k(u)+ \del_{jk}$, it follows that
\begin{equation}\label{sum}
\sum_{i=1}^s \ell_k(u_i)= \ell_k(u)+ (s-1)(n-\del_{jk}). 
\end{equation}
We next observe that if $w_i\in W^{(j)}$, then 
$$
\prod_{i=1}^m C_{w_i}
$$
is never a nonzero multiple of $C_\1$. Hence, after permuting the indices, for the elements $w_i$ as 
in the proposition, we obtain:
$$
w_{1},...,w_{m'}\in W^{(1)}, w_{m'+1},...,w_m\in W^{(2)}
$$ 
and for $m=m'+m''$, we have $1\le m', m''\le m-1$. Therefore,
\begin{equation}\label{products}
\prod_{i=1}^{m'} C_{w_i} =a' C_{w'}, \prod_{i=m'+1}^{m} C_{w_i} =a'' C_{w''},
\end{equation}
where $a', a''\ne 0$ in $\widehat{\k}$, and $w'\in W^{(1)}, w''\in W^{(2)}$. Moreover,
$$
C_{w'} C_{w''}=C_{\1}
$$
in $H_*(X, \widehat{\k})$. Therefore, by applying equations \eqref{sum} to the product decompositions 
\eqref{products}, we obtain 
$$
\sum_{k=1}^m \ell_k(w_i) = \ell_k(w')+\ell_k(w'') + mn -2n+1 -M
$$
where $M=m'\del_{1k}+ m''\del_{2k}\le m-1$. Since $\ell(w')+\ell(w'')=n$, it follows that
$$
\ell_k(w')+\ell_k(w'')= n-\del_{1k}-\del_{2k}=n-1. 
$$
Hence, we obtain
$$
\sum_{k=1}^m \ell_k(w_i) = (m-1)n -M\ge (m-1)(n-1). \qed 
$$

We are now ready to prove the theorem. We have to show that every 
$\ov{\la}=(\la_1,...,\la_m)\in \K_m$ satisfies the inequality
$$
 \sum_{i=1}^m w_i^{-1}(\la_i) \le_{\De^*} 0 
$$
for every $(w_1,...,w_m)\in \Si_{A,m}$. The latter is equivalent to two inequalities
$$
\sum_{i=1}^m \<\la_i, w_i(\zeta_k)\>=\sum_{i=1}^m \<w_i^{-1}(\la_i), \zeta_k\>\le 0, k=1, 2
$$
where $\zeta_k, k=1, 2$ are the vertices of the fundamental domain $\De_{sph}\subset S^1$ of $W$.

Suppose that this inequality fails for some $k$  and an $m$-tuple $(u_1,...,u_m)\in \Si_{A,m}$.
Since $\ov{\la}\in \K_m$, according to Theorem \ref{main}, 
there exists a semistable weighted configuration $\psi=(\mu_i\xi_i)$ in $X$ of the type $\ov{\la}$, 
so that the points $\xi_i$ belong to mutually antipodal spherical chambers 
$\De_{1},...,\De_{m}$ in $X$. Since 
$$
\prod_{i=1}^m C_{w_i}=C_{\1},
$$
for $r_i:=\ell_k(w_i)$, by combining Corollary \ref{C3} and Proposition \ref{ine}, we see that 
the intersection  
$$
\bigcap_{i=1}^m C_{r_i}(\De_{i}) 
$$
contains a vertex $\eta$ of type $\zeta_k$. Thus, as in the proof of Theorem \ref{main},
$$
slope_\psi(\eta)=-\sum_j \< \la_j, w_j(\zeta_k)\> <0.
$$
This contradicts semistability of $\psi$. \qed

\section{Triangle inequalities for associative commutative algebras}\label{generalize}

We now put the concept of stability inequalities in the more general context by associating a system of 
monoids $K_m(A)$ (generalizing the stability cones) to certain commutative and associative rings 
(which generalize the rings $A_t$). One advantage of this formalism is to eliminate dependence on the 
existence of the longest element $w_\circ\in W$ and getting more natural sets of inequalities. 
We also establish linear isomorphisms of cones $K_m(A_t)$ (defined ``cohomologically'') applying the above formalism to the algebras $A_t$  
and the Stability Cones $\K_{m+1}(Y)$ (defined ``homologically''). We conclude this section by showing 
that the system $WTI$ is irredundant. 

\medskip 
Let $\Lambda$ be a free abelian group (or a free module over an integral domain).

\begin{definition}
We say that a family of sub-monoids $K_m\subset \Lambda^{m+1}$, $m\ge 1$ is {\it coherent} if:

\noindent (1) The natural $S_m$-action on the first $m$-factors of $\Lambda^{m+1}$ preserves $K_m$;

\noindent (2) For any $(\lambda_1,\ldots,\lambda_m; \mu)\in \Lambda^{m+1}$ and $0<\ell<m$ the following are equivalent:

$\bullet$ $(\lambda_1,\ldots,\lambda_m;\mu)\in K_m$

$\bullet$ There exists $\mu'\in \Lambda$ such that
 $(\lambda_1,\ldots,\lambda_m;\mu')\in K_m$ and $(\mu',\lambda_{m+1},\ldots,\lambda_\ell;\mu)\in K_{\ell+1-m}$.

\end{definition}

Below we will interpret a coherent family of sub-monoids as above, as a commutative and associative (multivalued) operad.  

For any subsets $S\subset \Lambda^{m+1}=\Lambda^m\times \Lambda$, $S'\subset \Lambda^{\ell+1}=\Lambda\times \Lambda^\ell$ 
define the set $S'\circ S'\subset\Lambda^{m+k}=\Lambda^m\times  \Lambda^\ell$ to be the set of all
$(\lambda,\lambda')\in \Lambda^m\times  \Lambda^\ell$ such that there exists $\mu\in \Lambda$ such that $(\lambda,\mu)\in S$ and 
$(\mu,\lambda')\in S'$. In other words, if we regard $S, S'$ as correspondences $\Lambda^m\to \Lambda$ and $\Lambda\to \Lambda^\ell$, then 
$S\circ S'$ is their composition. The following is immediate:

\begin{lemma}
\label{le:circ coherence}
The second coherence condition is equivalent to that:
$$K_m\circ K_{\ell}= K_{m+\ell-1}$$
for all $m, \ell\ge 1$.

\end{lemma}

The following result is obvious.

\begin{lemma} If $K_m$, $m\ge 0$ is a coherent family, then each $K_m$, $m\ge 3$,  is the set of all 
$(\lambda_1,\ldots,\lambda_m;\mu)\in \Lambda^{m+1}$ such that there exist a sequence  $\mu_1,\ldots,\mu_m=\mu$ of elements in $\Lambda$ such that:
$(\lambda_1,\lambda_2;\mu_1)\in K_2$ and $(\mu_k,\lambda_{k+2};\mu_{k+1})\in K_2$ for $k=1,\ldots,m-1$.

\end{lemma}

We explain the naturality of the coherence condition below. To any submonoid $K_m\subset \Lambda^{m+1}$, 
$m\ge 1$ we associate an $m$-ary operation on subsets of $\Lambda$ as follows. For any subsets 
$S_1,S_2,\ldots,S_m\subset \Lambda$ define  $S_1\bullet S_2\bullet \cdot \ldots\bullet S_m\subset \Lambda_+$ 
to be image of the intersection $S_1\times \cdots\times S_m\times \Lambda\cap K_m$ 
under the projection to the $(m+1)$-st factor. 
In particular, if each $S_i=\{\lambda_i\}$ is  a single element set, then
$$\lambda_1\bullet \cdots \bullet \lambda_m=\{\mu\in \Lambda_+:(\lambda_1,\ldots,\lambda_m;\mu)\in K_m\} \ .$$
In general,
$$S_1\bullet \cdots\bullet S_m=\bigcup\limits_{(\lambda_1,\ldots,\lambda_m)\in S_1\times \cdots\times S_m} \lambda_1\bullet \cdots \bullet \lambda_m  \ .$$

\begin{lemma} If a family of submonoids $K_m\subset \Lambda^{m+1}$, $m\ge 1$ is  coherent, then  the above operations are:

\noindent (a)  commutative, i.e.,
$S_{\sigma(1)}\bullet \cdots \bullet S_{\sigma(m)}=S_1\bullet \cdots \bullet S_m$
for any permutation $\sigma$ of $\{1,\ldots,m\}$.

\noindent (b) Associative, i.e., 
$S_1\bullet \cdots \bullet S_k\bullet (S_{k+1}\bullet \cdots \bullet S_{\ell})\bullet S_{\ell+1}\bullet  \cdots \bullet S_m=S_1\bullet \cdots \bullet S_m$
for all $1\le k\le \ell \le m$ (i.e., informally speaking, these operations comprise a symmetric associative operad, see e.g., \cite{oper}).

\end{lemma}

\begin{proof} Part (a) is an obvious consequence of the first coherence condition.

We will now prove (b). Because of the already established commutativity, it suffices to verify the assertion for $k=0$. 
Also it suffices to proceed in the case when each $S_i=\{\lambda_i\}$ is an one-element set. That is, it suffices to prove that
$$(\lambda_1\bullet \cdots \bullet \lambda_\ell)\bullet \lambda_{\ell+1}\bullet \cdots \bullet \lambda_m=\lambda_1\bullet \cdots \bullet \lambda_m$$
The left hand side is the set of all $\mu\in \Lambda$ such that $(\mu',\lambda_{\ell+1},\ldots,\lambda_m;\mu)\in K_{m+1-\ell}$  
for some $\mu'\in \Lambda$ satisfying $(\lambda_1,\ldots,\lambda_\ell;\mu)\in K_\ell$.
By the second coherence condition, this is the set of all $\mu \in \Lambda$ such that $(\lambda_1,\ldots,\lambda_m;\mu)\in K_m$. 
But this set is exactly the right hand side of the above equation. This proves (b).

The lemma is proved. \qed 
\end{proof}

We now construct 
families of monoids associated with some associative commutative algebras.
Let  $\preceq$ be a partial order on $\Lambda$ compatible with the addition. 
This is equivalent to choosing a submonoid ${\mathcal M}$ ( the ``positive root cone'') 
such that $-{\mathcal M}\cap {\mathcal M}=\{0\}$,
so that $\lambda\preceq \mu$ if and only if  $\mu-\lambda\in {\mathcal M}$ 
(therefore, ${\mathcal M}=\{\lambda\in \Lambda:0\preceq \lambda\}$).

Let $A$ be commutative associative $\k$-algebra as in Section \ref{sec:BKfiltrations} with the basis 
labeled by a set $X\subset End(\Lambda)$ (i.e., the basis acts linearly on $\Lambda$). 
We define the structure coefficients $c_{x_1,\ldots,x_m}^y\in \k$ via
$$b_{x_1}\cdots  b_{x_m}= \sum_{y\in X} c_{x_1,\ldots,x_m}^y b_y
$$
for all $x_1,\ldots,x_n\in X$.

Given this data, we define:

$\bullet$ The dominant cone $\Lambda_+$ to be the set of all $\lambda\in \Lambda$ such that $x(\lambda)\preceq \lambda$ for all $x\in X$.

$\bullet$ For each $m\ge 0$ a subset $K_m(A)\subset \Lambda_+^{m+1}=\Lambda_+^m\times \Lambda_+$ to be the set of all  $(\la_1,...,\la_m;\mu)\in \Lambda_+^{m+1}$
such that
\begin{equation}\label{eq:cone}
y(\mu)\preceq x_1(\lambda_1)+\cdots +x_m(\lambda_m)
\end{equation}
for all $x_1,\ldots,x_m,y\in X$ such that  $c_{x_1,\ldots,x_m}^y\ne 0$ (with the convention that $K_0(A)=\Lambda_+$).

The following is immediate:

\begin{lemma}
\label{le:S_m action on K(A)}
The set $K_m(A)$ is a submonoid of $\Lambda^{m+1}$ invariant under the $S_m$-action on the first $m$ factors.
\end{lemma}

\begin{lemma}
\label{le:subcoherence}In the notation of Lemma \ref{le:circ coherence}  we have:
\begin{equation}
\label{eq:subcoherence}
K_m(A)\circ K_l(A)\subseteq K_{m+l-1}(A)
\end{equation}
for all $m,l\ge 1$.

\end{lemma}

\begin{proof} Indeed, let $(\lambda_1,\ldots,\lambda_{m+l-1};\mu)\in K_m(A)\circ K_l(A)$. 
This means that there exists $\mu_1\in \Lambda_+$ such that 
$(\lambda_1,\ldots,\lambda_m;\mu_1)\in K_m(A)$ and $(\mu_1,\lambda_{m+1},\ldots,\lambda_{m+l-1};\mu)\in K_l(A)$. Or, equivalently,
\begin{equation}\label{eq:cones}
y_1(\mu_1)\preceq x_1(\lambda_1)+\cdots +x_m(\lambda_m), y(\mu)\preceq y'_1(\mu_1) +x_{m+1}(\lambda_{m+1})+\cdots +x_{m+l-1}(\lambda_{m+l-1})
\end{equation}
for all $x_1,\ldots,x_{m+l-1},y_1,y\in X$ such that  $c_{x_1,\ldots,x_m}^{y_1}\ne 0$ and $c_{y'_1,x_{m+1},\ldots,x_{m+l-1}}^y\ne 0$. 
Now fix arbitrary $x_1,\ldots,x_{m+l-1},y\in X$ such that $c_{x_1,\ldots,x_{m+l-1}}^y\ne 0$. 
Due to associativity of multiplication in $A$ this implies existence of $y_1$ such that 
$c_{x_1,\ldots,m}^{y_1}\ne 0$ and $c_{y_1,x_{m+1},\ldots,x_{m+l-1}}^y\ne 0$. 
Therefore, we can take $y'_1=y_1$ in \eqref{eq:cones} and add the inequalities \eqref{eq:cones}. 
Hence, after canceling the term $y_1(\mu_1)$, we obtain
$$y(\mu)\preceq x_1(\lambda_1)+\cdots +x_{m+l-1}(\lambda_{m+l-1}) \ .$$
The latter inequality holds for all  $x_1,\ldots,x_{m+l-1},y\in X$ such that 
$c_{x_1,\ldots,x_{m+l-1}}^y\ne 0$ hence $(\lambda_1,\ldots,\lambda_{m+l-1};\mu)\in K_{m+l-1}(A)$. The lemma is proved. \qed 
\end{proof}

Thus, in view of Lemmas \ref{le:circ coherence}, \ref{le:S_m action on K(A)}, and \ref{le:subcoherence} the coherence  of 
$K_m(A)$, $m\ge 1$ depends entirely on whether or not the inclusion \eqref{eq:subcoherence} is an equality.

\begin{problem} Classify all commutative and associative algebras $A$ with basis labeled by 
$X\subset End(\Lambda)$ such that
\begin{equation}
\label{eq:supercoherence}
K_m(A)\circ K_l(A)\supseteq K_{m+l-1}(A)
\end{equation}
\end{problem}

We now specialize to the case associated with finite dihedral Weyl groups. Let $W=W_t$, where
$t=e^{\frac{\pi\sqrt{-1}}{n}}$, acting on the 2-dimensional real vector space $V$. We assume that 
$\R \otimes \La= V^*$; let ${\mathcal M}=\Delta^*\subset V^*$ be the dual cone to the positive (affine) Weyl chamber $\De\subset V$ of  $W$
(with respect to the simple roots $\al_1, \al_2$), i.e., $\De^*=\{\mu: \<\la, \mu\>\ge 0, \forall \la\in \De\}$). 
We take the based ring $A:=A_t$, with the basis $\{\sigma_w|w\in W_t\}$; accordingly, we take the based rings $B^{(i)}:=B_t^{(i)}, i=1, 2$.  
Thus, for $\zeta\in V$, $\la\in \R \otimes \La$, we have 
$$
\<\si_w(\la), \zeta\>=\<w^{-1}(\la), \zeta\>=\<\la, w(\zeta)\>. 
$$

Let $A_0, B_0^{(i)}$ be the associated graded algebras of $A, B^{(i)}$  with respect to the filtrations defined
by the concave function $\varphi, \varphi_i$ given by \eqref{eq:Belkale-Kumar convex function} as in
Theorem \ref{thm:filtrations}. Define 
$$
K_{m}(B)= K_{m}(B^{(1)}) \cap K_{m}(B^{(2)}), \quad K_{m}(B_0)= K_{m}(B^{(1)}_0) \cap K_{m}(B^{(2)}_0)
$$
Clearly, 
$$
K_{m}(A)\subset K_{m}(B) \subset K_{m}(B_0), \quad K_{m}(A)\subset K_{m}(A_0) \subset K_{m}(B_0)
$$

The following is the main result of the section. This is an analogue of the main result of \cite{bk} in 
the context of arbitrary finite dihedral groups. Recall that $\la^*=- w_\circ \la$. 

\begin{theorem} Assume that $t=e^{\frac{\pi\sqrt{-1}}{n}}$. Then for each $m\ge 2$ we have:
$$
K_{m}(B_0)=K_{m}(B)=K_m(A_0)= K_{m}(A).$$
Moreover, the above cones are isomorphic to the Stability Cone $\K_{m+1}(Y)$ for any thick spherical 
building $Y$ with the Weyl group $W$ via the 
linear map
$$
\Theta: (\la_1,...,\la_m; \mu)\mapsto (\mu_1=\la_1^*,...,\mu_m=\mu_m^*, \mu_{m+1}=\mu). 
$$ 
\end{theorem}
\proof Our goal is to relate the defining inequalities for the cone  $K_{m}(A)$ to Strong Triangle Inequalities; it will then follow that 
$$
K_{m}(B_0)=K_{m}(A). 
$$
Set $PD(w)=w_\circ w$ in $W$. 
Observe that for $u_1,...,u_m, v\in W$ and $\la_1,...,\la_m, \mu\in V^*$,
$$
\sum_{i=1}^m \si_{u_i}(\la_i) \ge_{\De^*} \si_v(\mu) \iff 
$$
$$
\sum_{i=1}^m u_i^{-1}(\la_i)\ge_{\De^*} - v^{-1} w_\circ \mu^*= - PD(v)^{-1} \mu^* \iff
$$
$$
\sum_{i=1}^m u_i^{-1}(\la_i) + PD(v)^{-1}\mu^* \ge_{\De^*} 0 \iff
$$
$$
\sum_{i=1}^{m+1} u_i^{-1}(\la_i) \ge_{\De^*} 0
$$
where $u_{m+1}:=PD(v)$ and $\la_{m+1}:=\mu^*$. Setting $w_i:=PD(u_i)$, $\mu_i:=\la_i^*$, we see that 
$$
\sum_{i=1}^{m+1} u_i^{-1}(\la_i) \ge_{\De^*} 0 \iff 
\sum_{i=1}^{m+1} w_i^{-1}(\mu_i) \le_{\De^*} 0. 
$$
Moreover, 
$$
c_{u_1,...,u_m}^v\ne 0 \hbox{~~in~~} A \iff 
$$
$$
c_{u_1,...,u_{m+1}}^{w_\circ}\ne 0  \hbox{~~in~~} A \iff 
$$
\begin{equation}\label{Xprod}
\prod_{i=1}^{m+1} C_{w_i}= a C_{\1}, a\ne 0  \hbox{~~in~~} H_*(Y, \widehat{\k}). 
\end{equation}
Recall that the system of inequalities
$$
\sum_{i=1}^{m+1} w_i^{-1}(\mu_i) \le_{\De^*} 0, \forall (w_1,...,w_{m+1}), \quad \hbox{so that~~} \eqref{Xprod} \hbox{~~holds} 
$$
is the system of Strong Triangle Inequalities. Therefore, the maps $u_i\mapsto PD(u_i)$, $i=1,...,m+1$, and 
$$
(\la_1,...,\la_m; \mu)\mapsto (\mu_1=\la_1^*,...,\mu_m=\mu_m^*, \mu_{m+1}=\mu)$$
determine a natural bijection between the set of defining inequalities for the cone $K_{m}(A)$ and the set of 
Strong Triangle Inequalities. Similarly, we obtain a bijection between the defining inequalities of $K_{m}(B_0)$ and 
the set of Weak Triangle Inequalities. However,  Strong Triangle Inequalities and Weak Triangle Inequalities determine the same cone, 
the Stability Cone $\K_{m+1}$, see Theorem \ref{equalcones}. Therefore, the map 
$$
\Theta: (\la_1,...,\la_m; \mu)\mapsto (\mu_1=\la_1^*,...,\mu_m=\mu_m^*, \mu_{m+1}=\mu)$$
determines linear isomorphisms of the cones
$$
K_{m}(A)\to \K_{m+1}, \quad K_{m}(B_0)\to  \K_{m+1}. 
$$
In particular, $K_{m}(A)=K_{m}(B)=K_{m}(B_0)$. Theorem follows. \qed 

\begin{cor}
$K_m(A)$ is invariant under $*: \la\to \la^*, \la\in \De$. 
\end{cor}
\proof Let $Y$ be a thick spherical building as above. Then $*$ extends to an isometry $*: Y\to Y$. Since isometries preserve (semi)stability 
condition, it follows that $\K_{m+1}=\K_{m+1}(Y)$ is invariant under $*$. Since $\Theta$ is $*$-equivariant, it follows that $K_m(A)$ is 
invariant under $*$ as well. \qed

\begin{cor}
For the algebra $A$ as above we have
$$
K_m(A)\circ K_l(A)= K_{m+l-1}(A). 
$$
\end{cor}
\proof Let ${\mathfrak Y}$ denote a thick Euclidean building modeled on $(\R^2, W)$. In view of the above theorem, we can interpret 
$K_k(A)$ as the set of $m+1$-tuples $(\la_1,...,\la_m; \mu)$ which are $\De$-valued side-lengths of ``disoriented'' 
geodesic $k+1$-gons $P=y_0\ldots y_{k}$ in ${\mathfrak Y}$, so that
$$
d_\De(y_{i-1}, y_{i})=\la_i, 1\le i\le k, \quad d_\De(y_{0}, y_{k})=\mu.
$$
(Note that the last side of $P$ has the orientation opposite to the rest.) 
For $k=m+l-1$, subdivide such a polygon by the diagonal $\ol{y_0 y_{l}}$ in two disoriented polygons
$$
P':=y_0 y_1\ldots y_{l}, \quad P''= y_0 y_l \ldots y_k.
$$
Then the $\De$-side lengths of these polygons are given by the tuples
$$
(\la_1,...,\la_{l}; \mu')\in K_l(A), \quad (\mu', \la_{l+1},..., \la_{k}; \mu'')\in K_m(A),
$$
where
$$
\mu'= d_\De(y_{0}, y_{l}), \quad \mu''=\mu.
$$
Hence, $K_{m+l-1}(A)\subset K_m(A)\circ K_l(A)$. \qed

\begin{thm}\label{irr}
The system of inequalities \eqref{wti} is irredundant.
\end{thm}
\proof The system of inequalities \eqref{wti} is nothing but the linear system
\begin{equation}\label{wti1}
\<w(\la_i - \la_j^*), \zeta_l\> \le \< \sum_{k\ne i, k\ne j} \la_k^*, \zeta_l\>, l=1, 2,
\quad w\in W.
\end{equation}

Fix regular vectors $\la_l\in \Delta, l=1,...,m, l\ne i, l\ne j$ 
(i.e., vectors from the interior of the Weyl chamber $\De$); then pick a vector $\la_j\in \De$ 
so that the distance from $\la_j$ to the boundary of $\De$ is at least $|\la_K^*|$, where
$$
\la_K^*:=\sum_{k\ne i, j} \la_k^*.
$$
Note that the vector $\la_K^*$ is again regular. 
Set 
$$
P=P_{\la_1,...,\la_{i-1}, \la_{i+1}, ... \la_m} = \la_j^*+ Hull(W\cdot \la_K^*)\subset \De.
$$
Here $Hull$ denotes the convex hull in $\R^2$. Then, for fixed  $\la_l, l\ne i$ as above, 
the solution set to the Weak Triangle Inequalities \eqref{wti} is exactly
the polygon $P$. Since $\la_K^*$ is regular, $P$ is a $2m$-gon.
Moreover, for each side of $P$ exactly one of the defining inequalities \eqref{wti1} is an equality.  \qed

\medskip
{\bf Beklale-Kumar inequalities.} In the context of complex algebraic reductive groups $G$, Belkale and 
Kumar \cite{bk} gave a certain description of the stability cone $\K_{m+1}$ using the rings $H_{BK}^*(G/P)$, where 
$P$ runs through the set of standard maximal parabolic subgroups of $G$, corresponding to the fundamental 
weights. In the context of rank 2 spherical buildings $X$, using our language, the system of 
Beklale-Kumar inequalities, imposed on vectors 
$$
(\la_1,...,\la_m; \mu)\in \De^{m+1},
$$
reads as follows: 
For every $(x_1,...,x_m;y)\in (W^{(k)})^{m+1}, k=1, 2$, so that
$c_{x_1,...,x_m}^y\ne 0$ in $gr(B^{(k)})$, we impose the inequality: 
$$
\sum_{i=1}^m \<\zeta_l, x_i(\la_i)\> \ge \<\zeta_l, y(\mu)\>.  
$$

We now observe that under the map $H_*(X_k, \widehat{\k})\to H_{BK}^*(X_l, \k):=gr(B^{(k)})$, 
determined by the inverse to the map \eqref{BKPR},  the ``infinities'' in $H_*(X_k, \widehat{\k})$ 
correspond to zeroes in $H_{BK}^*(X, \k)$. Accordingly, the structure constants equal to $1$ match structure
constants equal to $1$. Since the system $WTI$ is irredundant, we conclude that the system of 
Belkale-Kumar inequalities for $W=I_2(n)$, is also irredundant. 
Hence, Theorem \ref{irr} is an analogue (for $W=I_2(n)$) of a much deeper theorem by N.~Ressayre \cite{Re}, 
who proved irredundancy of Belkale-Kumar inequalities for arbitrary reductive groups.

Addresses:

Arkady Berenstein: Department of Mathematics, University of Oregon, Eugene, OR 97403,
USA. (arkadiy@uoregon.edu)

Michael Kapovich: Department of Mathematics, University of California, Davis, CA 95616,
USA. (kapovich@math.ucdavis.edu)

\end{document}